\documentclass[A4j,11pt]{article}
\usepackage{amsfonts,amsmath,amscd}
\usepackage{graphics}
\begin{document}

\title{{\bf A splitting theorem\\
for equifocal submanifolds\\
with non-flat section}}
\author{{\bf Naoyuki Koike}}
\date{}
\maketitle

\begin{abstract}
We first prove a certain kind of splitting theorem for 
an equifocal submanifold with non-flat section in a 
simply connected symmetric space of compact type, where an equifocal 
submanifold means a submanifold with parallel focal structure.  
By using the splitting theorem, we prove that each section of an equifocal 
submanifold with non-flat section in an irreducible simply connected 
symmetric space of compact type 
is isometric to a sphere or a real projective space.  
\end{abstract}

\vspace{0.5truecm}

\centerline{\bf 1. Introduction}

\vspace{0.5truecm}

A properly immersed complete submanifold $M$ in a simply connected 
symmetric space $G/K$ is called a {\it submanifold with parallel focal 
structure} if the following conditions hold:

\vspace{0.2truecm}

(PF-i) the restricted normal holonomy group of $M$ is trivial,

\vspace{0.2truecm}

(PF-ii) if $v$ is a parallel normal vector field on $M$ such that $v_{x_0}$ 
is a multiplicity $k$ focal normal of $M$ for some $x_0\in M$, then $v_x$ 
is a multiplicity $k$ focal normal of $M$ for all $x\in M$,

(PF-iii) for each $x\in M$, there exists a properly embedded complete 
connected submanifold through $x$ 
meeting all parallel submanifolds of $M$ orthogonally.  

\vspace{0.2truecm}

This notion was introduced by Ewert ([E2]).  
In [A], [AG] and [AT], 
this submanifold is simply called an {\it equifocal submanifold}.  
In this paper, we also shall use this name and 
assume that all equifocal submanifolds have trivial normal holonomy 
group.  The submanifold as in (PF-iii) ia called a {\it section of} 
$M$ {\it through} $x$, which is automatically totally geodesic.  
Note that Terng-Thorbergsson [TeTh] originally introduced the notion of 
an equifocal submanifold under the assumption that the sections is flat.  
The condition (PF-ii) is equivalent to the following condition:

\vspace{0.2truecm}

(PF-ii$'$) for each parallel unit normal vector field $v$ of $M$, the set of 
all focal radii along the geodesic $\gamma_{v_x}$ with $\dot{\gamma}_{v_x}(0)
=v_x$ is independent of the choice of $x\in M$.  

\vspace{0.2truecm}

\noindent
Note that, under the condition (PF-i), the condition (PF-iii) is equivalent to 
the following condition:

\vspace{0.2truecm}

(PF-iii$'$) $M$ has Lie triple systematic normal bundle (in the sense of 
[Koi1]).

\vspace{0.2truecm}

\noindent
In fact, (PF-iii)$\Rightarrow$(PF-iii$'$) is trivial and 
(PF-iii$'$)$\Rightarrow$(PF-iii) is shown as follows.  
If (PF-iii$'$) holds, then it is shown by Proposition 2.2 of [HLO] 
that $\exp^{\perp}(T^{\perp}_xM)$ meets all parallel submanifolds of $M$ 
orthogonally for each $x\in M$, where $\exp^{\perp}$ is the normal exponential 
map of $M$.  Also, it is clear that $\exp^{\perp}(T^{\perp}_xM)$ is properly 
embedded.  Thus (PF-iii) follows.  
An isometric action of a compact Lie group $H$ on a Riemannian manifold 
is said to be {\it polar} if there exists a properly embedded complete 
connected submanifold $\Sigma$ meeting every principal orbits of the 
$H$-action orthogonally.  The submanifold $\Sigma$ is called a {\it section} 
of the action.  If $\Sigma$ is flat, then the action is said to be 
{\it hyperpolar}.  
Principal orbits of polar actions are equifocal submanifolds 
and those of hyperpolar actions are equifocal ones with flat section.  
Conversely, homogeneous equifocal submanifolds (resp. 
homogeneous equifocal ones with flat section) in the symmetric spaces are 
catched as 
principal orbits of polar (resp. hyperpolar) actions on the spaces.  

In 1997, Heintze and Liu [HL] showed that an isoparametric submanifold in 
a Hilbert space is decomposed into a non-trivial (extrinsic) product of two 
such submanifolds if and only if the associated Coxeter group is 
decomposable.  In 1998, by using this splitting theorem of Heintze-Liu, 
Ewert [E1] showed that an equifocal submanifold with flat section 
in a simply connected 
symmetric space of compact type is decomposed into a non-trivial (extrinsic) 
product of two such submanifolds if and only if the associated Coxeter group 
is decomposable.  

In this paper, we first prove the following splitting theorem for an 
equifocal submanifold with non-flat section in a simply 
connected symmetric space of compact type.  

\vspace{0.3truecm}

{\bf Theorem A.} {\sl Let $M$ be an equifocal 
submanifold with non-flat section in a simply connected symmetric 
space $G/K$ of compact type and $\Sigma$ be a section of $M$.  
Then $M$ is decomposed into a non-trivial 
extrinsic product of two equifocal submanifolds if and only if 
the restricted holonomy group of 
{\rm(}the induced metric on{\rm)} $\Sigma$ is reducible.}

\vspace{0.3truecm}

Next we prove the following fact in terms of Theorem A.  

\vspace{0.3truecm}

{\bf Theorem B.} {\sl Let $M$ be an equifocal submanifold with non-flat 
section in an irreducible simply connected symmetric space $G/K$ of compact 
type.  Then each section of $M$ is isometric to a sphere or a real projective 
space.}

\vspace{1truecm}

\centerline{\bf Proof of Theorems A and B}

\vspace{0.5truecm}

In this section, we shall prove Theorem A.  
Without loss of generality, we may assume that 
$G$ is simply connected and $K$ is connected.  
Let $\pi:G\to G/K$ be the 
natural projection and $\phi:H^0([0,1],\mathfrak g)\to G$ be the parallel 
transport map for $G$, where $\mathfrak g$ is the Lie algebra of $G$ and 
$H^0([0,1],\mathfrak g)$ is the space of all $L^2$-integrable paths having 
$[0,1]$ as the domain.  Let $M^{\ast}:=\pi^{-1}(M)$ and $\widetilde M:=
(\pi\circ\phi)^{-1}(M)$.  Since $G$ is simply connected and $K$ is connected, 
$M^{\ast}$ and $\widetilde M$ are connected.  
Denote by $A$ (resp. $\widetilde A$) the shape tensor of $M$ (resp. 
$\widetilde M$) and by $\nabla^{\perp}$ (resp. $\widetilde{\nabla}^{\perp}$) 
the normal connection of $M$ (resp. $\widetilde M$).  
Let $\Sigma_x$ be the section of $M$ through $x\,(\in M)$.  
Assume that the restricted holonomy group of $\Sigma_x$ is reducible.  
Fix $x_0\in M$.  
We have the non-trivial orthogonal decomposition 
$T_x\Sigma_x=W_1\oplus W_2$, which is invariant with respect to 
the restricted holonomy group of $\Sigma_{x_0}$ at $x_0$.  
Since $M$ has trivial normal holonomy group, there exists the 
$\nabla^{\perp}$-parallel subbundle $D^N_i$ of the normal bundle $T^{\perp}M$ 
of $M$ with $(D^N_i)_{x_0}=W_i$ ($i=1,2$).  For each $x\in M$, it is easy to 
show that there exists an isometry $f$ of a neighborhood of $x_0$ in 
$\Sigma_{x_0}$ onto a neighborhood of $x$ in $\Sigma_x$ such that 
$f_{\ast x_0}$ coincides with the parallel translation (with respect to 
$\nabla^{\perp}$) along any curve in $M$ from $x_0$ to $x$.  From this fact, 
it follows that, for each $x\in M$, the orthogonal decomposition 
$T_x\Sigma_x=(D^N_1)_x\oplus(D^N_2)_x$ is invariant with respect to the 
restricted holonomy group of $\Sigma_x$ at $x$.  
Let $\widetilde D^N_i$ ($i=1,2$) be the subbundles of the normal 
bundle $T^{\perp}\widetilde M$ of $\widetilde M$ with 
$(\pi\circ\phi)_{\ast u}
((\widetilde D^N_i)_u)=(D^N_i)_{(\pi\circ\phi)(u)}$ ($u\in \widetilde M$) 
and $D^{N\ast}_i$ ($i=1,2$) be those of $T^{\perp}(M^{\ast})$ with 
$\pi_{\ast}((D^{N\ast}_i)_g)=(D^N_i)_{\pi(g)}$ ($g\in G$).  According to 
Lemma 1A.4 of [PoTh1], the focal set of $(M,x)$ consists of finitely many 
totally geodesic hypersurfaces in $\Sigma_x$.  Denote by $\mathfrak L_x$ the 
set of all focal hypersurfaces of $(M,x)$.  Let $\psi_x:\widehat{\Sigma}_x\to
\Sigma_x$ be the universal covering of $\Sigma_x$.  According to 
the de Rham's decomposition theorem, $\widehat{\Sigma}_x$ is isometric to 
the (non-trivial) Riemannian product 
$\widehat{\Sigma}^1_x\times\widehat{\Sigma}^2_x$, 
where $\widehat{\Sigma}^i_x$ ($i=1,2$) is the complete totally geodesic 
submanifold of $\widehat{\Sigma}_x$ through $\hat x\in\psi_x^{-1}(x)$ 
such that $(\psi_x)_{\ast\hat x}(T_{\hat x}\widehat{\Sigma}^i_x)
=(D^N_i)_x$.  By retaking the decomposition $T^{\perp}_{x_0}M=W_1\oplus W_2$ 
if necessary, we may assume that $\widehat{\Sigma}^1_x$ has no Euclidean part 
in the de Rham's decomposition for each $x\in M$.  
Let $\widehat{\mathfrak L}_x:=
\{\psi_x^{-1}(L)\,\vert\,L\in\mathfrak L_x\}$.  According to Corollary 3.6 
of [Kol2], elements of $\widehat{\mathfrak L}_x$ are either 
$L_1\times\widehat{\Sigma}^2_x$-type ($L_1\,:\,$ a totally geodesic 
hypersurface of 
$\widehat{\Sigma}^1_x$) or $\widehat{\Sigma}^1_x\times L_2$-type 
($L_2\,:\,$ a totally geodesic hypersurface of $\widehat{\Sigma}^2_x$), 
where we need the fact that $\widehat{\Sigma}^1_x$ has no Euclidean part.  
Denote by $\widehat{\mathfrak L}^1_x$ (resp. $\widehat{\mathfrak L}^2_x$) 
the set of all elements of $\widehat{\mathfrak L}_x$ of 
$L_1\times\widehat{\Sigma}^2_x$-type (resp. of 
$\widehat{\Sigma}^1_x\times L_2$-type) and set $\mathfrak L^i_x:=
\{L\in\mathfrak L_x\,\vert\,\psi_x^{-1}(L)\in\widehat{\mathfrak L}^i_x\}$ 
($i=1,2$).  Let $V':=
\overline{{\rm Span}(\displaystyle{\mathop{\cup}_{u\in\widetilde M}
T^{\perp}_u\widetilde M)}}$, $V_i:=
\overline{{\rm Span}(\displaystyle{\mathop{\cup}_{u\in\widetilde M}
(\widetilde D^N_i)_u)}}$ ($i=1,2$) and $V_0:=(V')^{\perp}$.  Also, let 
$(\widetilde D^T_0)_u:=\displaystyle{\mathop{\cap}_{v\in T^{\perp}_u
\widetilde M}{\rm Ker}\,\widetilde A_v}$, 
$(\widetilde D^T_1)_u
:=\displaystyle{\left(\mathop{\cap}_{v\in (\widetilde D^N_2)_u}{\rm Ker}\,
\widetilde A_v\right)\ominus(\widetilde D^T_0)_u}$ and 
$(\widetilde D^T_2)_u
:=\displaystyle{\left(\mathop{\cap}_{v\in (\widetilde D^N_1)_u}{\rm Ker}\,
\widetilde A_v\right)\ominus(\widetilde D^T_0)_u}$, where $u\in\widetilde M$.  
Without loss of generality, we may assume that $\widetilde M$ 
includes the zero element 
$\hat 0$ of $H^0([0,1],\mathfrak g)$, where we note that $\hat 0$ is the 
constant path at the zero element $0$ of $\mathfrak g$.  
Let $\widetilde M':=\widetilde M\cap V'$.  
First we prepare the following fact.  

\vspace{0.5truecm}

{\bf Proposition 2.1.} 
{\sl We have $\widetilde M=\widetilde M'\times V_0\subset 
V'\times V_0=H^0([0,1],\mathfrak g)$.}

\vspace{0.5truecm}

{\it Proof.} First we shall show $V_0\subset(\widetilde D^T_0)_u$ for each 
$u\in\widetilde M$, where we regard $(\widetilde D^T_0)_u(\subset 
T_uH^0([0,1],\mathfrak g))$ as a subspace of $H^0([0,1],\mathfrak g)$ under 
the identification of $T_uH^0([0,1],\mathfrak g)$ with $H^0([0,1],
\mathfrak g)$.  
From the definition of $V_0$, we have $V_0\subset T_u\widetilde M$ for each 
$u\in\widetilde M$.  Let $(\widetilde D^T_0)_u^{\perp}$ be the orthogonal 
complement of $(\widetilde D^T_0)_u$ in $T_u\widetilde M$.  Clearly we have 
$(\widetilde D^T_0)_u^{\perp}=\displaystyle{\sum_{v\in T^{\perp}_u
\widetilde M}\left(\mathop{\oplus}_{\lambda\in{\rm Spec}\,\widetilde A_v
\setminus\{0\}}{\rm Ker}(\widetilde A_v-\lambda\,{\rm id})\right)}$, where 
${\rm Spec}\,\widetilde A_v$ is the spectrum of $\widetilde A_v$.  Let 
$X\in{\rm Ker}(\widetilde A_v-\lambda\,{\rm id})$ ($v\in T^{\perp}_u
\widetilde M,\,\lambda\in{\rm Spec}\,\widetilde A_v\setminus\{0\}$).  
Let $J_X$ be the strongly Jacobi field along the normal geodesic 
$\gamma_v$ with $\gamma_v'(0)=v$ satisfying $J_X(0)=X$ (hence 
$J_X'(0)=-A_vX$).  
Let $\alpha:(-\varepsilon,\varepsilon)\to M$ be a curve in $M$ with 
$\alpha'(0)=X$ and $\widetilde v$ be the parallel normal vector field along 
$\alpha$ with $\widetilde v_0=v$.  Define a map 
$\delta:(-\varepsilon,\varepsilon)\times[0,\infty)\to H^0([0,1],\mathfrak g)$ 
by $\delta(t,s):=\gamma_{\widetilde v_t}(s)$, where $\gamma_{\widetilde v_t}$ 
is the normal geodesic in $H^0([0,1],\mathfrak g)$ with 
$\gamma'_{\widetilde v_t}(0)=\widetilde v_t$.  Then we have 
$\delta_{\ast}(\frac{\partial}{\partial t}\vert_{t=0})=J_X$.  
Since $\delta(t,0)-\delta(t,\frac{1}{\lambda})\in T^{\perp}_{\alpha(t)}
\widetilde M\subset V'$ for each $t\in(-\varepsilon,\varepsilon)$, we have 
$\delta_{\ast}(\frac{\partial}{\partial t}\vert_{t=s=0})-\delta_{\ast}
(\frac{\partial}{\partial t}\vert_{t=0,s=\frac{1}{\lambda}})\in V'$.  
On the other hand, we have 
$\delta_{\ast}(\frac{\partial}{\partial t}\vert_{t=s=0})=X$ and 
$\delta_{\ast}(\frac{\partial}{\partial t}\vert_{t=0,s=\frac{1}{\lambda}})=0$.  Hence we have $X\in V'$.  From the arbitrariness of $X$, it follows that 
${\rm Ker}(\widetilde A_v-\lambda\,{\rm id})\subset V'$.  Furthermore, it 
follows from the arbitrarinesses of $\lambda$ and $v$ that 
$(\widetilde D^T_0)_u^{\perp}\subset V'$, that is, $V_0\subset
(\widetilde D_0^T)_u$.  Since $V_0\subset(\widetilde D^T_0)_u\subset 
T_u\widetilde M$ for any $u\in \widetilde M$, we have 
$\widetilde M=\displaystyle{\mathop{\cup}_{u\in \widetilde M'}(u+V_0)
=\widetilde M'\times V_0\subset V'\times V_0}$.\hspace{0.1truecm}q.e.d.

\vspace{0.5truecm}

\centerline{
\unitlength 0.1in
\begin{picture}( 40.2800, 27.2500)( -0.6200,-29.1600)
%
\special{pn 8}%
\special{pa 926 2916}%
\special{pa 3326 2916}%
\special{fp}%
%
\special{pn 8}%
\special{pa 3326 2916}%
\special{pa 3966 1948}%
\special{fp}%
\special{pa 3966 1948}%
\special{pa 1566 1948}%
\special{fp}%
\special{pa 1566 1948}%
\special{pa 926 2916}%
\special{fp}%
%
\special{pn 8}%
\special{ar 2846 3076 1440 1440  3.2537955 4.0476145}%
%
\special{pn 8}%
\special{pa 1406 2916}%
\special{pa 1406 1156}%
\special{fp}%
%
\special{pn 8}%
\special{ar 2830 1324 1440 1440  3.2537955 4.0476145}%
%
\special{pn 8}%
\special{pa 1950 1948}%
\special{pa 1950 196}%
\special{fp}%
%
\special{pn 8}%
\special{pa 1726 1156}%
\special{pa 3166 1476}%
\special{fp}%
\special{pa 1566 1396}%
\special{pa 3166 1476}%
\special{fp}%
%
\special{pn 8}%
\special{pa 1726 1156}%
\special{pa 1502 1396}%
\special{fp}%
\special{sh 1}%
\special{pa 1502 1396}%
\special{pa 1562 1362}%
\special{pa 1538 1358}%
\special{pa 1534 1334}%
\special{pa 1502 1396}%
\special{fp}%
\special{pa 2142 1252}%
\special{pa 2022 1420}%
\special{fp}%
\special{sh 1}%
\special{pa 2022 1420}%
\special{pa 2078 1378}%
\special{pa 2054 1378}%
\special{pa 2044 1354}%
\special{pa 2022 1420}%
\special{fp}%
%
\special{pn 8}%
\special{pa 2542 1340}%
\special{pa 2486 1444}%
\special{fp}%
\special{sh 1}%
\special{pa 2486 1444}%
\special{pa 2536 1396}%
\special{pa 2512 1398}%
\special{pa 2500 1376}%
\special{pa 2486 1444}%
\special{fp}%
%
\special{pn 8}%
\special{pa 1590 2380}%
\special{pa 1590 620}%
\special{fp}%
%
\special{pn 8}%
\special{pa 1726 1316}%
\special{pa 1726 996}%
\special{fp}%
\special{pa 1878 900}%
\special{pa 1878 628}%
\special{fp}%
\special{pa 1726 844}%
\special{pa 1726 524}%
\special{fp}%
%
\special{pn 20}%
\special{sh 1}%
\special{ar 3166 1476 10 10 0  6.28318530717959E+0000}%
\special{sh 1}%
\special{ar 1590 2372 10 10 0  6.28318530717959E+0000}%
\special{sh 1}%
\special{ar 1590 2372 10 10 0  6.28318530717959E+0000}%
%
\special{pn 8}%
\special{pa 2206 516}%
\special{pa 1726 588}%
\special{fp}%
\special{sh 1}%
\special{pa 1726 588}%
\special{pa 1796 598}%
\special{pa 1780 580}%
\special{pa 1790 558}%
\special{pa 1726 588}%
\special{fp}%
\special{pa 2206 516}%
\special{pa 1878 732}%
\special{fp}%
\special{sh 1}%
\special{pa 1878 732}%
\special{pa 1946 712}%
\special{pa 1924 704}%
\special{pa 1924 680}%
\special{pa 1878 732}%
\special{fp}%
\special{pa 2206 516}%
\special{pa 1726 1124}%
\special{fp}%
\special{sh 1}%
\special{pa 1726 1124}%
\special{pa 1784 1084}%
\special{pa 1760 1082}%
\special{pa 1752 1060}%
\special{pa 1726 1124}%
\special{fp}%
%
\special{pn 8}%
\special{pa 1294 748}%
\special{pa 1590 828}%
\special{fp}%
\special{sh 1}%
\special{pa 1590 828}%
\special{pa 1532 792}%
\special{pa 1540 814}%
\special{pa 1520 830}%
\special{pa 1590 828}%
\special{fp}%
\put(17.6600,-23.0800){\makebox(0,0)[rt]{$\hat 0$}}%
\put(37.1800,-20.1200){\makebox(0,0)[rt]{$V'$}}%
\put(12.9400,-6.2800){\makebox(0,0)[rt]{$V_0$}}%
\put(23.9800,-4.0400){\makebox(0,0)[rt]{$\widetilde D^T_0$}}%
\put(15.5800,-11.8800){\makebox(0,0)[rt]{$X$}}%
\put(22.3800,-14.6000){\makebox(0,0)[rt]{$J_X$}}%
%
\special{pn 8}%
\special{pa 1758 2676}%
\special{pa 1510 2556}%
\special{fp}%
\special{sh 1}%
\special{pa 1510 2556}%
\special{pa 1562 2604}%
\special{pa 1558 2580}%
\special{pa 1580 2568}%
\special{pa 1510 2556}%
\special{fp}%
\put(19.8200,-26.2000){\makebox(0,0)[rt]{${\widetilde M}'$}}%
\put(11.9800,-14.2000){\makebox(0,0)[rt]{$\widetilde M$}}%
%
\special{pn 8}%
\special{ar 2054 1636 536 616  3.5395008 4.0450204}%
\end{picture}%
\hspace{1.5truecm}
}

\vspace{0.5truecm}

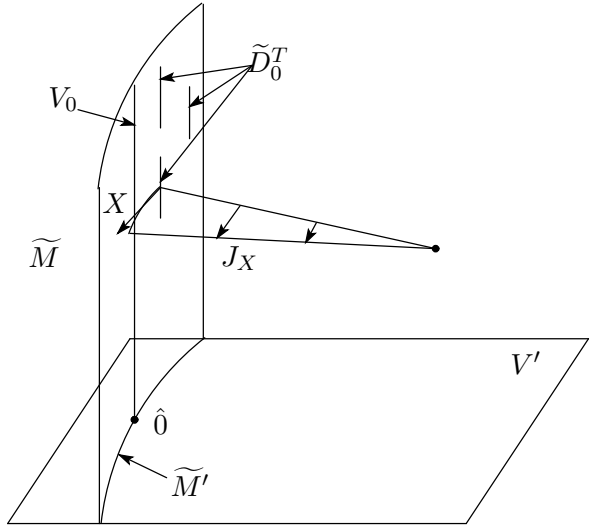
\centerline{{\bf Fig. 1.}}

\noindent
Define distributions $D^T_0,\,D^T_1$ and $D^T_2$ on $M$ by 
$$\begin{array}{l}
\displaystyle{(D^T_0)_x:=\left(\mathop{\cap}_{v\in T^{\perp}_xM}{\rm Ker}
\,A_v\right)\cap g_{\ast}\left(\mathfrak c_{g_{\ast}^{-1}T_xM}(g_{\ast}^{-1}
T^{\perp}_xM)\right),}\\
\displaystyle{(D^T_1)_x:=\left(\left(\mathop{\cap}_{v\in (D^N_2)_x}
{\rm Ker}\,A_v\right)\cap g_{\ast}\left(\mathfrak c_{g_{\ast}^{-1}T_xM}
(g_{\ast}^{-1}(D^N_2)_x)\right)\right)\ominus(D^T_0)_x,}\\
\displaystyle{(D^T_2)_x:=\left(\left(\mathop{\cap}_{v\in(D^N_1)_x}
{\rm Ker}\,A_v\right)\cap g_{\ast}\left(\mathfrak c_{g_{\ast}^{-1}T_xM}(g_{\ast}^{-1}
(D^N_1)_x)\right)\right)\ominus(D^T_0)_x,}
\end{array}$$
for each $x=gK\in M$, where 
$\mathfrak c_{\ast}(\sharp)$ is the centralizer of $\sharp$ in $\ast$.  
Take an arbitrary $v\in T^{\perp}_{eK}M$.  
Set $\mathfrak p:=T_{eK}(G/K)$.  Denote by ${\rm Spec}\,R(v)$ the 
spectrum of $R(v):=R(\cdot,v)v$, where $R$ is the curvature tensor of $G/K$.  
For $\mu\in{\rm Spec}R(v)$, we set $\mathfrak p_{\mu}^v:=
{\rm Ker}(R(v)-\mu{\rm id})$, $\mathfrak f_{\mu}^v:={\rm ad}(v)
\mathfrak p_{\mu}^v$ ($\mu\in{\rm Spec}R(v)\setminus\{0\}$) and 
$\mathfrak f_0^v:=({\rm Ker}\,{\rm ad}(v))\cap\mathfrak f$.  
Note that 
$$T_{eK}M=\mathfrak p^v_0\cap T_{eK}M
+\sum_{\mu\in{\rm Spec}\,R(v)\setminus\{0\}}
(\mathfrak p^v_{\mu}\cap T_{eK}M)\leqno{(2.2)}$$
and 
$$T^{\perp}_{eK}M=\mathfrak p^v_0\cap T^{\perp}_{eK}M+
\sum_{\mu\in{\rm Spec}\,R(v)\setminus\{0\}}
(\mathfrak p^v_{\mu}\cap T^{\perp}_{eK}M)\leqno{(2.3)}$$
because $M$ is equifocal and hence it has Lie triple 
systematic normal bundle.  
Denote by $\mathfrak f$ the Lie algebra of $K$.  
For each $X\in\mathfrak g$, we define loop vectors ${\it l}^{co}_{X,k}$ and 
${\it l}^{si}_{x,k}$ ($k\in{\bf N}$) by 
${\it l}^{co}_{X,k}(t):=X\cos(2k\pi t)$ and 
${\it l}^{si}_{X,k}(t):=X\sin(2k\pi t)$.  
For $X\in\mathfrak p_{\mu}^v$ 
$(\mu\in{\rm Spec}\,R(v)\setminus\{0\}$), 
we set $X_{\mathfrak f}:=\frac{1}{\sqrt{\mu}}{\rm ad}(v)(X)$.  
For $X\in\mathfrak p_{\mu}^v$, $Y\in\mathfrak p_0^v\oplus\mathfrak f_0^v$ and 
$k\in{\bf Z}$, we define loop vectors ${\it l}^i_{v,X,k},\,
{\it l}^i_{v,X_{\mathfrak f},k}$ and ${\it l}^i_{v,Y,k}\in 
H^0([0,1],\mathfrak g)$ ($i=1,2$) by 
$$\begin{array}{l}
\displaystyle{{\it l}^1_{v,X,k}(t)={\it l}^1_{v,X_{\mathfrak f},k}(t)
={\it l}^{co}_{X,k}(t)-{\it l}^{si}_{X_{\mathfrak f},k}(t),}\\
\displaystyle{{\it l}^2_{v,X,k}(t)={\it l}^2_{v,X_{\mathfrak f},k}(t)=
={\it l}^{si}_{X,k}(t)+{\it l}^{co}_{X_{\mathfrak f},k}(t).}
\end{array}$$
For a general $Z\in\mathfrak g$, we define loop vectors 
${\it l}^i_{v,Z,k}\in H^0([0,1],\mathfrak g)$ ($i=1,2,\,k\in{\bf Z}$) by 
$${\it l}^i_{v,Z,k}:={\it l}^i_{v,Z_0,k}
+\sum_{\mu\in{\rm Spec}\,R(v)\setminus\{0\}}
\left({\it l}^i_{v,Z_{\mathfrak p,\mu},k}+
{\it l}^i_{v,Z_{\mathfrak f,\mu},k}\right),$$
where $Z=Z_0+\sum\limits_{\mu\in{\rm Spec}\,R(v)\setminus\{0\}}
(Z_{\mathfrak p,\mu}+Z_{\mathfrak f,\mu})$ 
($Z_0\in\mathfrak p_0^v,\,Z_{\mathfrak p,\mu}\in
\mathfrak p_{\mu}^v,\,Z_{\mathfrak f,\mu}\in
\mathfrak f_{\mu}^v$).  Denote by $\widehat{\ast}$ the constant path 
at $\ast\in\mathfrak g$.  Note that $\widehat{\ast}$ is the horizontal lift 
of $\ast\,(\in\mathfrak g=T_eG)$ to $\widehat 0$.  
Then, according to Propositions 3.1 and 3.2 of [Koi2] and those proofs, 
we have the following relations.  

\vspace{0.5truecm}

{\bf Lemma 2.2.} {\sl Let $X\in T_{eK}M\cap\mathfrak p_{\mu}^v$.  
Then we have 
$$\begin{array}{l}
\displaystyle{\tilde A_{\hat v}{\it l}^1_{v,X,k}
=\frac{\sqrt{\mu}}{2k\pi}(\hat X-{\it l}^1_{v,X,k}),}\\
\displaystyle{\tilde A_{\hat v}{\it l}^2_{v,X,k}
=\frac{\sqrt{\mu}}{2k\pi}(\hat X_{\mathfrak f}-{\it l}^2_{v,X,k}),}\\
\displaystyle{\tilde A_{\hat v}\hat X=\widehat{A_vX}-\frac{\sqrt{\mu}}{2}
\hat X_{\mathfrak f}+\frac{\sqrt{\mu}}{2\pi}\sum_{k\in{\bf Z}\setminus\{0\}}
\frac{1}{k}{\it l}^1_{v,X,k},}\\
\displaystyle{\tilde A_{\hat v}\hat X_{\mathfrak f}
=-\frac{\sqrt{\mu}}{2}\hat X+\frac{\sqrt{\mu}}{2\pi}
\sum_{k\in{\bf Z}\setminus\{0\}}\frac{1}{k}{\it l}^2_{v,X,k}}
\end{array}$$
and 
$$\widetilde{\nabla}^{\perp}_{{\it l}^1_{v,X,k}}\widetilde v^L
=\widetilde{\nabla}^{\perp}_{{\it l}^2_{v,X,k}}\widetilde v^L
=\widetilde{\nabla}^{\perp}_{\hat X}\widetilde v^L
=\widetilde{\nabla}^{\perp}_{\hat X_{\mathfrak f}}\widetilde v^L
=0,$$
where $k\in{\bf Z}\setminus\{0\}$ and 
$\widetilde v^L$ is the horizontal lift of 
a parallel normal vector field $\widetilde v$ with $\widetilde v_0=v$ 
along an arbitrary curve $\alpha$ in $M$ with $\dot{\alpha}(0)=X$.}

\vspace{0.5truecm}

{\bf Lemma 2.3.} {\sl Let $w\in T^{\perp}_{eK}M\cap\mathfrak p_{\alpha}^v$.  
Then we have 
$$\begin{array}{l}
\displaystyle{\tilde A_{\hat v}{\it l}^1_{v,w,k}
=-\frac{\sqrt{\mu}}{2k\pi}{\it l}^1_{v,w,k},}\\
\displaystyle{\tilde A_{\hat v}{\it l}^2_{v,w,k}
=\frac{\sqrt{\mu}}{2k\pi}(\hat w_{\mathfrak f}-{\it l}^2_{v,w,k}),}\\
\displaystyle{\tilde A_{\hat v}\hat w_{\mathfrak f}=\frac{\sqrt{\mu}}{2\pi}
\sum_{k\in{\bf Z}\setminus\{0\}}\frac{1}{k}
{\it l}^2_{v,w,k},}
\end{array}$$
and 
$$\widetilde{\nabla}^{\perp}_{{\it l}^1_{v,w,k}}\widetilde v^L
=-\frac{\sqrt{\mu}}{2k\pi}\hat w,\,\,
\widetilde{\nabla}^{\perp}_{{\it l}^2_{v,w,k}}\widetilde v^L=0,\,\,
\widetilde{\nabla}^{\perp}_{\hat w_{\mathfrak f}}\widetilde v^L
=\frac{\sqrt{\mu}}{2}\hat w,$$
where $k\in{\bf Z}\setminus\{0\}$ and 
$\widetilde v^L$ is as in Lemma 2.2.}

\vspace{0.5truecm}

{\bf Lemma 2.4.} {\sl Let $X\in \mathfrak p^v_0$ and 
$Y\in\mathfrak f^v_0$.  Then we have 
$$\tilde A_{\hat v}{\it l}^{co}_{X,k}
=\tilde A_{\hat v}{\it l}^{si}_{X,k}
=\tilde A_{\hat v}{\it l}^{co}_{Y,k}
=\tilde A_{\hat v}{\it l}^{si}_{Y,k}
=\tilde A_{\hat v}\hat Y=0,$$
$$\widetilde{\nabla}^{\perp}_{{\it l}^{co}_{X,k}}\widetilde v^L
=\widetilde{\nabla}^{\perp}_{{\it l}^{si}_{X,k}}\widetilde v^L
=\widetilde{\nabla}^{\perp}_{{\it l}^{co}_{Y,k}}\widetilde v^L
=\widetilde{\nabla}^{\perp}_{{\it l}^{si}_{Y,k}}\widetilde v^L
=\widetilde{\nabla}^{\perp}_{\hat Y}\widetilde v^L=0$$
and
$$\tilde A_{\hat v}\hat X=\widehat{A_vX},\,\,\,\,
\widetilde{\nabla}^{\perp}_{\widehat X}\widetilde v^L=0\,\,\,\,
({\rm when}\,\,X\in\mathfrak p^v_0\cap T_{eK}M),$$
where $i=1,2,\,k\in{\bf N}$ and $\widetilde v^L$ is as in Lemma 2.2.}

\vspace{0.5truecm}

From Lemmas $2.2\sim2.4$, we can show the following relations.  

\vspace{0.5truecm}

{\bf Lemma 2.5.} {\sl At $\hat 0\in\widetilde M$, $(\widetilde D^T_0)_{\hat0}$ 
is equal to 
$$\begin{array}{l}
\displaystyle{{\rm Span}\{\widehat X\,\vert\,
X\in(D^T_0)_{eK}\}\oplus{\rm Span}\{\widehat{\eta}\,\vert\,\eta\in \mathfrak c
_{\mathfrak f}(T^{\perp}_{eK}M)\}}\\
\hspace{0.6truecm}\displaystyle{\oplus{\rm Span}\{{\it l}^{co}_{Z,k}\,\vert\,
Z\in \mathfrak c_{\mathfrak g}(T^{\perp}_{eK}M),\,k\in{\bf N}\setminus
\{0\}\}}\\
\hspace{0.6truecm}\displaystyle{\oplus{\rm Span}\{{\it l}^{si}_{Z,k}\,\vert\,
Z\in \mathfrak c_{\mathfrak g}(T^{\perp}_{eK}M),\,k\in{\bf N}\setminus
\{0\}\}}
\end{array}$$
and $(\widetilde D^T_{j_1})_{\hat0}$ is equal to 
$$\begin{array}{l}
\displaystyle{{\rm Span}\{\widehat X\,\vert\,
X\in(D^T_{j_1})_{eK}\}\oplus{\rm Span}\{\widehat{\eta}\,\vert\,\eta\in 
\mathfrak c_{\mathfrak f}((D^N_{j_2})_{eK})\ominus \mathfrak c_{\mathfrak f}
(T^{\perp}_{eK}M)\}}\\
\hspace{0.6truecm}\displaystyle{\oplus{\rm Span}\{{\it l}^{co}_{Z,k}\,\vert\,
Z\in \mathfrak c_{\mathfrak g}((D^N_{j_2})_{eK})\ominus \mathfrak c
_{\mathfrak g}(T^{\perp}_{eK}M),\,k\in{\bf N}\setminus\{0\}\}}\\
\hspace{0.6truecm}\displaystyle{\oplus{\rm Span}\{{\it l}^{si}_{Z,k}\,\vert\,
Z\in \mathfrak c_{\mathfrak g}((D^N_{j_2})_{eK})\ominus \mathfrak c
_{\mathfrak g}(T^{\perp}_{eK}M),\,k\in{\bf N}\setminus\{0\}\},}
\end{array}$$
where $(j_1,j_2)=(1,2)$ or $(2,1)$.}

\vspace{0.5truecm}

{\it Proof.} 
According to Lemmas $2.2\sim2.4$, we have 
$$\begin{array}{l}
\displaystyle{{\rm Ker}\,\widetilde A_{\hat v}
={\rm Span}
\{\widehat X\,\vert\,X\in{\rm Ker}\,A_v\cap \mathfrak p^v_0\}\oplus
{\rm Span}\{\widehat{\eta}\,\vert\,\eta\in\mathfrak f_0^v\}}\\
\hspace{1.6truecm}\displaystyle{\oplus{\rm Span}\{{\it l}^{co}_{Z,k}\,\vert\,
Z\in{\rm Ker}\,{\rm ad}(v),\,k\in{\bf N}\setminus\{0\}\}}\\
\hspace{1.6truecm}\displaystyle{\oplus{\rm Span}\{{\it l}^{si}_{Z,k}\,\vert\,
Z\in{\rm Ker}\,{\rm ad}(v),\,k\in{\bf N}\setminus\{0\}\}.}
\end{array}$$
Hence we have the desired relations.  \hspace{4.4truecm}q.e.d.

\vspace{0.5truecm}

From Lemmas $2.2\sim2.4$, we have the following lemma.  

\vspace{0.5truecm}

\noindent
{\bf Lemma 2.6.} {\sl Assume that $v\in D^N_i$.  Let $\widetilde v^L$ be as in 
Lemma 2.2.  Then the statements {\rm(i)} and {\rm (ii)} hold.  

{\rm (i)} For each $X\in T\widetilde M$, we have $\widetilde{\nabla}^{\perp}_X
\widetilde v^L\in\widetilde D^N_i$.  

{\rm (ii)} For each $Y\in \widetilde D^T_j\oplus\widetilde D^T_0$ ($j\not=i$), 
we have $\widetilde{\nabla}_Y\widetilde v^L=0$.}

\vspace{0.5truecm}

\noindent
{\it Proof.} Without loss of generality, we may assume that the base point 
of $X$ is $\hat 0$.  
First we shall show the statement (i).  According to $(2.2),\,(2.3)$ and 
Lemmas $2.2\sim2.4$, we have only to show $\widetilde{\nabla}^{\perp}_X
\widetilde v^L\in\widetilde D^N_i$ in case of $X={\it l}^1_{v,w,k}$ or 
$\widehat w_{\mathfrak f}$ ($w\in T^{\perp}_{eK}M\cap\mathfrak p_{\mu}^v$).  
Suppose that $w\in D^N_j$ ($j\not=i$).  Then we have 
$\mu=0$ because the sectional curvature of ${\rm Span}\{v,w\}$ is equal 
to $0$.  This contradicts $\mu\not=0$.  
Hence we have $w\in D^N_i$.  Therefore, it follows from Lemma 2.3 that 
$\widetilde{\nabla}^{\perp}_X\widetilde v^L\in{\rm Span}\{\widehat w\}
\subset\widetilde D^N_i$.  Thus the statement (i) is shown.  
Next we shall show the statement (ii).  
From (i), we have $\widetilde{\nabla}^{\perp}_Y
\widetilde v^L=0$.  Also, from the definitions of $\widetilde D^T_j$ and 
$\widetilde D^T_0$, we have $\widetilde A_{\hat v}Y=0$.  Hence, we obtain 
$\widetilde{\nabla}_Y\widetilde v^L=0$.  
\hspace{5.5truecm}q.e.d.

\vspace{0.5truecm}

By using (ii) of Lemma 2.6, we prove the following lemma.  

\vspace{0.5truecm}

{\bf Lemma 2.7.} {\sl For each $u\in\widetilde M$, the tangent space 
$T_u\widetilde M$ is orthogonally decomposed as $T_u\widetilde M=
(\widetilde D^T_1)_u\oplus(\widetilde D^T_2)_u\oplus(\widetilde D^T_0)_u$.}

\vspace{0.5truecm}

{\it Proof.} Take unit vectors $v_i$ belonging to $(\widetilde D^N_i)_u$ 
($i=1,2$).  According to (i) of Lemma 2.6, we have 
$\widetilde R^{\perp}(X,Y)v_1\in(\widetilde D^N_1)_u$ 
for any $X,Y\in T_u\widetilde M$, where 
$\widetilde R^{\perp}$ is the curvature tensor of the normal connection of 
$\widetilde M$.  Hence, it follows from the Ricci equation that 
$[\widetilde A_{v_1},\widetilde A_{v_2}]=0$.  Therefore, we have 
$$T_u\widetilde M=\mathop{\oplus}_{\lambda\in{\rm Spec}\,\widetilde A_{v_1}}
\mathop{\oplus}_{\mu\in{\rm Spec}\,\widetilde A_{v_2}}\left({\rm Ker}
(\widetilde A_{v_1}-\lambda\,{\rm id})\cap{\rm Ker}(\widetilde A_{v_2}
-\mu\,{\rm id})\right), \leqno{(2.1)}$$
where ${\rm Spec}\,\widetilde A_{v_i}$ ($i=1,2$) is the spectrum of 
$\widetilde A_{v_i}$.  Set $\widetilde{\mathfrak L}_u^i:=\{(\pi\circ\phi)
\vert_{T^{\perp}_u\widetilde M})^{-1}(L)\,\vert\,L\in{\mathfrak L}^i_{\pi(u)}
\}$ ($i=1,2$).  The family $\widetilde {\mathfrak L}^1_u\cup
\widetilde{\mathfrak L}^2_u$ gives the family of all focal hypersurfaces of 
$\widetilde M$ at $u$.  Let $\lambda\in{\rm Spec}\,\widetilde A_{v_1}
\setminus\{0\}$ and $\mu\in{\rm Spec}\,\widetilde A_{v_2}\setminus\{0\}$.  
We shall show ${\rm Ker}(\widetilde A_{v_1}-\lambda\,{\rm id})\cap
{\rm Ker}(\widetilde A_{v_2}-\mu\,{\rm id})=\{0\}$.  Suppose that ${\rm Ker}
(\widetilde A_{v_1}-\lambda\,{\rm id})\cap{\rm Ker}(\widetilde A_{v_2}-
\mu\,{\rm id})\not=\{0\}$.  Take $X(\not=0)\,\in{\rm Ker}(\widetilde A_{v_1}
-\lambda\,{\rm id})\cap{\rm Ker}(\widetilde A_{v_2}-\mu\,{\rm id})$.  
The point $u+\frac{1}{\lambda}v_1$ and $u+\frac{1}{\mu}v_2$ are focal points 
along the normal geodesics $\gamma_{v_1}$ and $\gamma_{v_2}$, respectively.  
Hence there exist $L_1\in\widetilde{\mathfrak L}^1_u$ with 
$u+\frac{1}{\lambda}v_1\in L_1$ and $L_2\in\widetilde{\mathfrak L}^2_u$ 
with $u+\frac{1}{\mu}v_2\in L_2$.  Let $w_{\theta}:=\cos\,\theta\cdot v_1+
\frac{\lambda}{\mu}\sin\,\theta\cdot v_2$ ($0\leq\theta\leq\frac{\pi}{2}$).  
Since $A_{w_{\theta}}X=\lambda(\sin\,\theta+\cos\,\theta)X$, the point 
$u+\frac{1}{\lambda(\sin\,\theta+\cos\,\theta)}w_{\theta}$ is a focal point 
along $\gamma_{w_{\theta}}$ for each $\theta\in[0,\frac{\pi}{2}]$.  Define 
a curve $c:[0,\frac{\pi}{2}]\to H^0([0,1],\mathfrak g)$ by $c(\theta):=u+
\frac{1}{\lambda(\sin\,\theta+\cos\,\theta)}w_{\theta}$ ($\theta\in I$), 
which is smooth and regular.  For each $\theta\in[0,\frac{\pi}{2}]$, we have 
$c(\theta)\in\displaystyle{\mathop{\cup}_{L\in\widetilde{\mathfrak L}^1_u
\cup\widetilde{\mathfrak L}^2_u}(L\cap{\rm Span}\{v_1,v_2\})}$.  
For simplicity, we set $F:=\displaystyle{\mathop{\cup}_{L\in\widetilde
{\mathfrak L}^1_u\cup\widetilde{\mathfrak L}^2_u}(L\cap
{\rm Span}\{v_1,v_2\})}$.  Since $F$ is a family of affine lines in 
${\rm Span}\{v_1,v_2\}$ which are parallel to ${\rm Span}\{v_1\}$ or 
${\rm Span}\{v_2\}$ and $c$ is a regular curve in $F$, $c$ lies in the 
only affine line belonging to $F$.  It is clear that the affine lines 
$L_1\cap{\rm Span}\{v_1,v_2\}$ and $L_2\cap{\rm Span}\{v_1,v_2\}$ are mutually 
distinct.  These facts contradict $c(0)\in L_1$ and 
$c(\frac{\pi}{2})\in L_2$ (see Fig. 2).  Therefore we have 
${\rm Ker}(\widetilde A_{v_1}-\lambda\,{\rm id})\cap{\rm Ker}
(\widetilde A_{v_2}-\mu\,{\rm id})=\{0\}$.  This fact together with $(2.1)$ 
deduces 
$\displaystyle{\mathop{\oplus}_{\lambda\in{\rm Spec}\,\widetilde A_{v_1}
\setminus\{0\}}{\rm Ker}(\widetilde A_{v_1}-\lambda\,{\rm id})\subset{\rm Ker}
\,\widetilde A_{v_2}}$.  
From the arbitrariness of $v_2$, we have 
$\displaystyle{\mathop{\oplus}_{\lambda\in
{\rm Spec}\,\widetilde A_{v_1}\setminus\{0\}}{\rm Ker}(\widetilde A_{v_1}-
\lambda\,{\rm id})\subset(\widetilde D^T_0)_u\oplus(\widetilde D^T_1)_u}$.  
That is, the orthogonal complement 
$((\widetilde D^T_0)_u\oplus(\widetilde D^T_1)_u)^{\perp}$ of 
$(\widetilde D^T_0)_u\oplus(\widetilde D^T_1)_u$ is contained in 
${\rm Ker}\,\widetilde A_{v_1}$.  From the arbitrariness of $v_1$, we have 
$((\widetilde D^T_0)_u\oplus(\widetilde D^T_1)_u)^{\perp}\subset
(\widetilde D^T_0)_u\oplus(\widetilde D^T_2)_u$, which implies 
$((\widetilde D^T_0)_u\oplus(\widetilde D^T_1)_u)^{\perp}\subset
(\widetilde D^T_2)_u$.  On the other hand, we have 
$((\widetilde D^T_0)_u\oplus(\widetilde D^T_1)_u)\cap
(\widetilde D^T_2)_u=\{0\}$.  Hence we have $T_u\widetilde M=
(\widetilde D^T_0)_u\oplus(\widetilde D^T_1)_u\oplus(\widetilde D^T_2)_u$ and 
$((\widetilde D^T_0)_u\oplus(\widetilde D^T_1)_u)^{\perp}
=(\widetilde D^T_2)_u$.  After all we have $T_u\widetilde M
=(\widetilde D^T_0)_u\oplus(\widetilde D^T_1)_u\oplus(\widetilde D^T_2)_u$ 
(orthogonal direct sum).  
\begin{flushright}q.e.d.\end{flushright}

\vspace{0.5truecm}

\centerline{
\unitlength 0.1in
\begin{picture}( 28.1700, 17.4600)( 15.3300,-24.2600)
%
\special{pn 8}%
\special{pa 2334 968}%
\special{pa 2334 2228}%
\special{fp}%
\special{pa 2334 2228}%
\special{pa 4134 2228}%
\special{fp}%
\special{pa 4134 2228}%
\special{pa 4134 968}%
\special{fp}%
\special{pa 4134 968}%
\special{pa 2334 968}%
\special{fp}%
%
\special{pn 8}%
\special{pa 2694 2048}%
\special{pa 3064 2048}%
\special{fp}%
\special{sh 1}%
\special{pa 3064 2048}%
\special{pa 2996 2028}%
\special{pa 3010 2048}%
\special{pa 2996 2068}%
\special{pa 3064 2048}%
\special{fp}%
\special{pa 2694 2048}%
\special{pa 2694 1688}%
\special{fp}%
\special{sh 1}%
\special{pa 2694 1688}%
\special{pa 2674 1756}%
\special{pa 2694 1742}%
\special{pa 2714 1756}%
\special{pa 2694 1688}%
\special{fp}%
%
\special{pn 8}%
\special{pa 3774 2228}%
\special{pa 3774 968}%
\special{fp}%
\special{pa 4134 1328}%
\special{pa 2334 1328}%
\special{fp}%
%
\special{pn 8}%
\special{pa 3064 2048}%
\special{pa 3774 2048}%
\special{da 0.070}%
\special{pa 2694 1688}%
\special{pa 2694 1328}%
\special{da 0.070}%
%
\special{pn 8}%
\special{pa 3774 2048}%
\special{pa 2694 1328}%
\special{fp}%
%
\special{pn 8}%
\special{pa 4350 1508}%
\special{pa 3774 1680}%
\special{fp}%
\special{sh 1}%
\special{pa 3774 1680}%
\special{pa 3844 1680}%
\special{pa 3826 1664}%
\special{pa 3832 1642}%
\special{pa 3774 1680}%
\special{fp}%
\special{pa 3568 852}%
\special{pa 3352 1328}%
\special{fp}%
\special{sh 1}%
\special{pa 3352 1328}%
\special{pa 3398 1276}%
\special{pa 3374 1280}%
\special{pa 3360 1260}%
\special{pa 3352 1328}%
\special{fp}%
%
\special{pn 8}%
\special{pa 3414 1662}%
\special{pa 3334 1752}%
\special{fp}%
\special{sh 1}%
\special{pa 3334 1752}%
\special{pa 3392 1716}%
\special{pa 3370 1712}%
\special{pa 3364 1688}%
\special{pa 3334 1752}%
\special{fp}%
\put(31.7100,-18.7700){\makebox(0,0)[rt]{$v_1$}}%
\put(26.5800,-17.0600){\makebox(0,0)[rt]{$v_2$}}%
\put(54.2100,-13.3700){\makebox(0,0)[rt]{$L_1\cap{\rm Span}\{v_1,v_2\}$}}%
\put(42.8700,-6.8000){\makebox(0,0)[rt]{$L_2\cap{\rm Span}\{v_1,v_2\}$}}%
\put(35.9400,-15.1700){\makebox(0,0)[rt]{$c(\theta)$}}%
\put(36.0300,-24.2600){\makebox(0,0)[rt]{${\rm Span}\{v_1,v_2\}$}}%
\end{picture}%
\hspace{1truecm}
}

\vspace{0.5truecm}

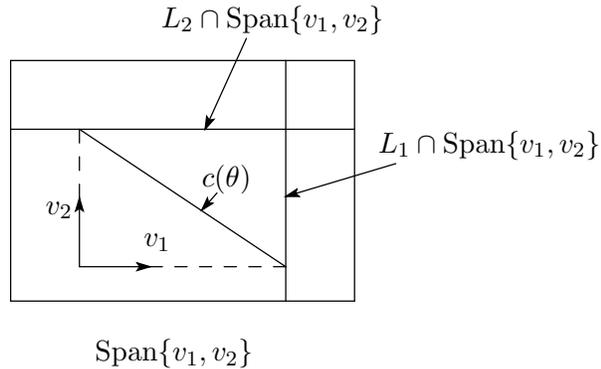
\centerline{{\bf Fig. 2.}}

\vspace{0.5truecm}

Next we prepare the following lemma.  

\vspace{0.5truecm}

{\bf Lemma 2.8.} {\sl {\rm (i)} The distributions 
$\widetilde D^T_i\oplus\widetilde D^T_0$ {\rm(}$i=1,2${\rm)} 
are totally geodesic.

{\rm (ii)} The distributions $\widetilde D^T_i$ {\rm(}$i=1,2${\rm)} 
are totally geodesic.}

\vspace{0.5truecm}

{\it Proof.} For simplicity, set $\widetilde D^T_{i0}:=\widetilde D^T_i
\oplus\widetilde D^T_0$ ($i=1,2$).  Denote by $\widetilde h$ (resp. 
$\widetilde h_{10}$) the second fundamental form of $\widetilde M$ 
(resp. $\widetilde D^T_{10}$), by $\widetilde A^{10}$ the shape tensor of 
$\widetilde D^T_{10}$, by $\widetilde{\nabla}$ (resp. $\nabla^{\widetilde M}$) 
the Levi-Civita connection of $H^0([0,1],\mathfrak g)$ (resp. $\widetilde M$) 
and by $\nabla^{\perp_2}$ the normal connection of $\widetilde D^T_2$.  
Also, denote by $\overline{\nabla}$ the connection of the bundle 
$T^{\ast}\widetilde M\otimes T^{\ast}\widetilde M\otimes T^{\perp}
\widetilde M$ induced from $\nabla^{\widetilde M}$ and 
$\widetilde{\nabla}^{\perp}$.  
Let $X,Y\in(\widetilde D^T_{10})_u$ and $Z\in(\widetilde D^T_2)_u$.  Let 
$\widetilde X$ (resp. $\widetilde Y$) be a section of $\widetilde D^T_{10}$ 
with $\widetilde X_u=X$ (resp. $\widetilde Y_u=Y$) and $\widetilde Z$ be a 
section of $\widetilde D^T_2$ with $\widetilde Z_u=Z$.  For any 
$v_1\in(\widetilde D^N_1)_u$, we have 
$\langle\widetilde h(Y,Z),v_1\rangle=\langle\widetilde A_{v_1}Z,Y\rangle=0$ 
because of $(\widetilde D^T_2)_u\subset{\rm Ker}\,\widetilde A_{v_1}$.  
Also, for any $v_2\in(\widetilde D^N_2)_u$, we have 
$\langle\widetilde h(Y,Z),v_2\rangle=\langle\widetilde A_{v_2}Y,Z\rangle=0$ 
because of $(\widetilde D^T_{10})_u\subset{\rm Ker}\,\widetilde A_{v_2}$.  
Hence we have $\widetilde h(Y,Z)=0$.  From the arbitrarinesses of $Y,Z$ and 
$u$, we have $\widetilde h(\widetilde D^T_{10},\widetilde D^T_2)=0$.  
Also, we can show $\widetilde h(\widetilde D^T_{10},\widetilde D^T_{10})
\subset\widetilde D^N_1$ and $\widetilde h(\widetilde D^T_2,\widetilde D^T_2)
\subset\widetilde D^N_2$.  
Let $X,Y,Z,\widetilde Y$ and $\widetilde Z$ be as above.  It follows from 
$\widetilde h(\widetilde D^T_{10},\widetilde D^T_2)=0$ that 
$$\begin{array}{l}
\displaystyle{(\overline{\nabla}_X\widetilde h)(Z,Y)=
\widetilde{\nabla}^{\perp}_X
(\widetilde h(\widetilde Z,\widetilde Y))-\widetilde h(
\nabla^{\widetilde M}_X\widetilde Z,Y)
-\widetilde h(Z,\nabla^{\widetilde M}_X\widetilde Y)}\\
\hspace{2.3truecm}\displaystyle{=\widetilde h(A^{10}_ZX,Y)-\widetilde h(Z,
h_{10}(X,Y))}\\
\hspace{2.3truecm}\displaystyle{\equiv-\widetilde h(Z,h_{10}(X,Y))\qquad
({\rm mod}\,\,(\widetilde D^N_1)_u).}
\end{array}\leqno{(2.2)}$$
Also, it follows from $\widetilde h(\widetilde D^T_{10},\widetilde D^T_{10})
\subset\widetilde D^N_1$ and Lemma 2.6 that 
$$\begin{array}{l}
(\overline{\nabla}_Z\widetilde h)(X,Y)=\widetilde{\nabla}^{\perp}_Z
(\widetilde h(\widetilde X,\widetilde Y))-\widetilde h(\nabla^{\perp_2}_Z
\widetilde X,Y)-\widetilde h(X,\nabla^{\perp_2}_Z\widetilde Y)\\
\hspace{2.3truecm}\equiv 0\,\,
({\rm mod}\,(\widetilde D^N_1)_u).
\end{array}
\leqno{(2.3)}$$
By $(2.2),\,(2.3)$ and the Codazzi equation, we have 
$\widetilde h(Z,h_{10}(X,Y))\in(\widetilde D^N_1)_u$.  
On the other hand, it follows from $\widetilde h(\widetilde D^T_2,
\widetilde D^T_2)\subset\widetilde D^N_2$ that 
$\widetilde h(Z,h_{10}(X,Y))\in(\widetilde D^N_2)_u$.  Hence we have 
$\widetilde h(Z,h_{10}(X,Y))=0$.  
According to the proof of Lemma 2.7, we have 
$$(\widetilde D^T_2)_u=
\displaystyle{\mathop{\oplus}_{v_2\in(\widetilde D^N_2)_u}
\mathop{\oplus}_{\mu\in{\rm Spec}\widetilde A_{v_2}\setminus\{0\}}{\rm Ker}
(\widetilde A_{v_2}-\mu\,{\rm id})}.$$
If $Z\in{\rm Ker}(\widetilde A_{v_2}-\mu\,{\rm id})$ ($\mu\in{\rm Spec}
\widetilde A_{v_2}\setminus\{0\}$), then we have 
$$\langle\widetilde h(Z,h_{10}(X,Y)),v_2\rangle=\langle\widetilde A_{v_2}Z,
h_{10}(X,Y)\rangle=\mu\langle h_{10}(X,Y),Z\rangle=0,$$
that is, 
$\langle h_{10}(X,Y),Z\rangle=0$.  
From the arbitrariness of $Z\in(\widetilde D^T_2)_u$, it follows that 
$h_{10}(X,Y)=0$.  
From the arbitrarinesses of $X$ and $Y$, it follows that $h_{10}=0$, that is, 
$\widetilde D^T_{10}$ is totally geodesic.  Similarly, we can show that 
$\widetilde D^T_{20}$ is totally geodesic.  By the similar discussion, 
we can show the statement (ii).  
\hspace{8.4truecm}q.e.d.

\vspace{0.5truecm}

By using Lemmas 2.6$\sim$2.8, we show the following fact.  

\vspace{0.5truecm}

{\bf Lemma 2.9.} {\sl We have $V'=V_1\oplus V_2$ $($orthogonal direct sum$)$.}

\vspace{0.5truecm}

{\it Proof.} Clearly we have $V'=V_1+V_2$.  We have only to show 
$V_1\perp V_2$.  Take arbitrary $u_1,u_2\in\widetilde M$ and arbitrary 
$v_j\in(\widetilde D^N_j)_{u_i}$ ($(i,j)=(1,2),(2,1)$).  Define a subset 
$U(u_1)$ of $H^0([0,1],\mathfrak g)$ by $U(u_1):=$\newline
$\displaystyle{\mathop{\cup}_{u\in L^{\widetilde D^T_{10}}_{u_1}}
L^{\widetilde D^T_{20}}_u}$, where $L^{\widetilde D^T_{10}}_{u_1}$ 
(resp. $L^{\widetilde D^T_{20}}_u$) 
is the leaf of $\widetilde D^T_{10}$ (resp. 
$\widetilde D^T_{20}$) through $u_1$ (resp. $u$).  Since $\widetilde M$ is 
complete, $\widetilde D^T_{10}$ is totally geodesic by Lemma 2.8 and 
$\widetilde D^T_2$ is the orthogonal complementary distribution of 
$\widetilde D^T_{10}$ by Lemma 2.7, $\widetilde D^T_2$ is an Ehresmann 
connection for the 
foliation consisting of integral manifolds of $\widetilde D^T_{10}$ (see 
[BH]).  Note that the discussions in [BH] are valid in the infinite 
dimensional case.  From the infinite dimensional version of the discussion 
in [BH], it follows that $U(u_1)=\widetilde M$.  Therefore we have 
$L^{\widetilde D^T_{10}}_{u_1}\cap L^{\widetilde D^T_{20}}_{u_2}\not=
\emptyset$.  Take $u_3\in L^{\widetilde D^T_{10}}_{u_1}\cap 
L^{\widetilde D^T_{20}}_{u_2}$ and curves 
$\alpha_i:[0,1]\to L^{\widetilde D^T_{i0}}_{u_i}$ ($i=1,2$) with 
$\alpha_i(0)=u_i$ and $\alpha_i(1)=u_3$.  According to (ii) of Lemma 2.6, 
we have 
$P^{\widetilde{\nabla}}_{\alpha_i}(v_j)\in(\widetilde D^N_j)_{u_3}$ 
($i=1,2$), where $P^{\widetilde{\nabla}}_{\alpha_i}$ is the parallel 
translation along $\alpha_i$ with respect to $\widetilde{\nabla}$.  Hence 
we obtain $\langle v_1,v_2\rangle=0$.  Therefore, it follows from the 
arbitrarinesses of $v_1$ and $v_2$ that $V_1\perp V_2$.  
\hspace{8.7truecm}q.e.d.

\vspace{0.5truecm}

Fix $x_0\in M$.  According to Lemma 2.15 and Proposition 2.16 of [E2], 
the focal set of 
$(M,x_0)$ consists of finitely many totally geodesic hypersurfaces in the 
section $\Sigma_{x_0}$ through $x_0$.  Let $\mathfrak L_{x_0}$ be the family 
of all the focal hypersurfaces.  The focal hypersurfaces divide $\Sigma_{x_0}$ 
into some open cells.  Denote by $\triangle$ the component containing 
$0\in T^{\perp}_{x_0}M$ of the inverse image by 
$\exp^{\perp}_{x_0}$ of the open cell containing $x_0$.  Define a map 
$f:M\times\triangle\to G/K$ by $f(x,v):=\exp^{\perp}_x(\widetilde v_x)$ 
($(x,v)\in M\times\triangle$), where $\widetilde v$ is the parallel normal 
vector field of $M$ with $\widetilde v_{x_0}=v$.  Let $U:=f(M\times\triangle)$, which is an open dense subset of $G/K$ consisting of non-focal 
points of $M$.  For each $v\in\triangle$, denote by $M_v$ 
the parallel submanifold $\eta_{\widetilde v}(M)$ of $M$, where 
$\eta_{\widetilde v}$ is the end-point map for $\widetilde v$, 
that is, $\eta_{\widetilde v}(x)=f(x,v)$ ($x\in M$).  
Let $E^N_i$ ($i=1,2$) be the distribution on $U$ such that $E^N_i\vert_M=D^N_i,
\,\,E^N_i\vert_{\Sigma_x}$ is a parallel distribution on $\Sigma_x$ for each 
$x\in M$ and that $E^N_i\vert_{M_v}$ is a normal parallel subbundle of 
$T^{\perp}M_v$ for each $v\in\triangle$.  
Denote by $(D^T_i)^v$ ($i=0,1,2$) the distributions on 
$M_v$ corresponding to the distributions $D^T_i$ ($i=0,1,2$) on $M$.  
It is shown that $(D^T_i)^v=(\eta_{\widetilde v})_{\ast}(D^T_i)$.  
For each $i\in\{0,1,2\}$, the distributions $(D^T_i)^v$'s ($v\in\triangle$) 
give a distribution on $U$.  
Denote by $E^T_i$ ($i=0,1,2$) this distribution on $U$.  Set 
$E_i:=E_i^T\oplus E_i^N$ and $E_{i0}:=E_i^T\oplus E_i^N\oplus E_0^T$ 
($i=1,2$).  Let $\widetilde U:=(\pi\circ\phi)^{-1}(U)$, which is an open 
dense subset of $H^0([0,1],\mathfrak g)$.  For each $v\in\triangle$, denote 
by $\widetilde M_v$ the submanifold $\eta_{\widetilde v^L}(\widetilde M)$, 
where $\eta_{\widetilde v^L}$ is the end-point map for the horizontal lift 
$\widetilde v^L$ of $\widetilde v$.  
Note that $\eta_{\widetilde v^L}(\widetilde M)$ is not a parallel submanifold 
of $\widetilde M$ because $\widetilde v^L$ is not parallel with respect to 
the normal connection of $\widetilde M$.  Let $\widetilde E^N_i$ ($i=1,2$) be the horizontal lift of $E^N_i$ to $\widetilde U$.  
Denote by $(\widetilde D^T_i)^v$ the distributions on $\widetilde M_v$ 
corresponding to the distributions $\widetilde D^T_i$ ($i=0,1,2$) on 
$\widetilde M$.  For each $i\in\{0,1,2\}$, the distributions 
$(\widetilde D^T_i)^v$'s ($v\in\triangle$) give a distribution on 
$\widetilde U$.  Denote by 
$\widetilde E^T_i$ ($i=0,1,2$) this distribution.  Set 
$\widetilde E_i:=\widetilde E^T_i\oplus\widetilde E^N_i$ and 
$\widetilde E_{i0}:=\widetilde E_i^T\oplus\widetilde E_i^N\oplus
\widetilde E^T_0$ ($i=1,2$).  
By using Lemmas 2.5 and 2.8, we show the following lemma.  

\vspace{0.5truecm}

{\bf Lemma 2.10.} {\sl {\rm (i)} The distributions $\widetilde E_{i0}$ 
{\rm(}$i=1,2${\rm)} are totally geodesic.

{\rm (ii)} The distributions $\widetilde E_i$ {\rm(}$i=1,2${\rm)} 
are totally geodesic.}

\vspace{0.5truecm}

{\it Proof.} For each $X\in\Gamma(TM)$, we define $\overline X\in\Gamma(TU)$ 
by $\overline X_{f(x,v)}:=(\eta_{\widetilde v})_{\ast x}(X_x)$ 
($(x,v)\in M\times\triangle$), where $\eta_{\widetilde v}$ is as above.  
Also, for each $w\in\triangle$, we define $\overline w\in\Gamma(TU)$ by 
$\overline w_{f(x,v)}:=P^{\Sigma_x}_{\gamma_{\widetilde v_x}}
(\widetilde w_x)$ ($(x,v)\in M\times\triangle$), where $\widetilde w$ is 
the parallel normal vector field of $M$ with $\widetilde w_{x_0}=w$ and 
$P^{\Sigma_x}_{\gamma_{\widetilde v_x}}$ is the parallel translation along 
the geodesic $\gamma_{\widetilde v_x}:[0,1]\to\Sigma_x$ with 
$\gamma_{\widetilde v_x}'(0)=\widetilde v_x$ with respect to the Levi-Civita 
connection of $\Sigma_x$.  Note that $P^{\Sigma_x}_{\gamma_{\widetilde v_x}}$ 
coincides with the parallel translation along $\gamma_{\widetilde v_x}$ with 
respect to the Levi-Civita connection of $G/K$ because $\Sigma_x$ is totally 
geodesic.  Without loss of generality, we may assume $x_0=eK$.  We suffice 
to show that $\widetilde E_{i0}$ ($i=1,2$) and $\widetilde E_i$ ($i=1,2$) 
have $\hat 0$ as a geodesic point.  
Easily we can show that 
if $X\in\Gamma(D^T_i)$ (resp. $w\in\triangle\cap(D^N_j)_{eK}$), then 
$\overline X\in\Gamma(E_i^T)$ (resp. $\overline w\in\Gamma(E^N_j)$), where 
$i=0,1,2$ and $j=1,2$.  We shall show that $\widetilde E_{10}$ has $\hat 0$ 
as a geodesic point.  
From Lemma 2.5, we have 
$$\begin{array}{l}
\displaystyle{(\widetilde E_{10})_{\hat0}={\rm Span}\{\widehat X\,\vert\,
X\in(D^T_{10})_{eK}\}\oplus{\rm Span}\{\widehat{\eta}\,\vert\,\eta\in
\mathfrak c_{\mathfrak f}((D^N_2)_{eK})\}}\\
\hspace{1.5truecm}\displaystyle{\oplus{\rm Span}\{{\it l}^{co}_{Z,k}\,\vert\,
Z\in\mathfrak c_{\mathfrak g}((D^N_2)_{eK}),\,k\in{\bf N}\setminus\{0\}\}}\\
\hspace{1.5truecm}\displaystyle{\oplus{\rm Span}\{{\it l}^{si}_{Z,k}\,\vert\,
Z\in\mathfrak c_{\mathfrak g}((D^N_2)_{eK}),\,k\in{\bf N}\setminus\{0\}\}}\\
\hspace{1.5truecm}\displaystyle{\oplus{\rm Span}\{\widehat w\,\vert\,w\in
(D^N_1)_{eK}\}.}
\end{array}\leqno{(2.4)}$$
Denote by $\widetilde h_{10}$ the second fundamental form of 
$\widetilde E_{10}$.  First we show $\widetilde h_{10}((\widetilde D^N_1)_{\hat0},(\widetilde D^N_1)_{\hat 0})=0$.  Let $w_1,w_2\in(D^N_1)_{eK}$.  Denote 
by $\nabla,\,\nabla^{\ast}$ and $\widetilde{\nabla}$ the Levi-Civita 
connection of $G/K, G$ and $H^0([0,1],\mathfrak g)$.  Denote by $(\cdot)^L$ 
(resp. $(\cdot)^{\ast}$) the horizontal lift of $(\cdot)$ to $H^0([0,1],
\mathfrak g)$ (resp. $G$).  According to Lemmas 2.2 and 2.3 in [Koi2], we have 
$$\widetilde{\nabla}_{\widehat w_1}\overline w_2^L=(\nabla^{\ast}_{w_1}
\overline w_2^{\ast})_{\hat0}^L-t[w_1,w_2]+\frac12[w_1,w_2]_{\hat 0}^L
=(\nabla_{w_1}\overline w_2)_{\hat0}^L-t[w_1,w_2],$$
where $t[w_1,w_2]$ is the $H^0$-path in $\mathfrak g$ assigning 
$t[w_1,w_2]$ to each $t\in[0,1]$.  Since $E^N_1$ is totally geodesic, 
we have $\nabla_{w_1}\overline w_2\in(D^N_1)_{eK}$ and hence 
$(\nabla_{w_1}\overline w_2)_{\hat 0}^L\in(\widetilde E_{10})_{\hat0}$ by 
$(2.4)$.  Also, we have $[w_1,w_2]\in\mathfrak c_{\mathfrak f}((D^N_2)_{eK})$ 
and hence $t[w_1,w_2]\in(\widetilde E_{10})_{\hat0}$ by $(2.4)$.  Therefore, 
we have $\widetilde{\nabla}_{\widehat w_1}\overline w_2^L\in
(\widetilde E_{10})_{\hat 0}$, that is, $\widetilde h_{10}(\widehat w_1,
\widehat w_2)=0$.  Thus we have 
$$\widetilde h_{10}((\widetilde D^N_1)_{\hat0},(\widetilde D^N_1)_{\hat0})=0.
\leqno{(2.5)}$$
Set $\widetilde E^T_{10}:=\widetilde E^T_1\oplus\widetilde E^T_0$.  Next we 
show that $\widetilde h_{10}((\widetilde E^T_{10})_{\hat0},
(\widetilde E^T_{10})_{\hat0})=0$.  Let $\widetilde X,\,\widetilde Y
\in\Gamma(\widetilde E^T_{10})$.  
For each $w\in(D^N_2)_{eK}$, we have 
$$\langle\widetilde h(\widetilde X_{\hat0},\widetilde Y_{\hat0}),\widehat w
\rangle=\langle\widetilde A_{\widehat w}\widetilde X_{\hat0},
\widetilde Y_{\hat 0}\rangle=0$$
from the definition of $\widetilde E^T_{10}$.  Hence we have 
$\widetilde h(\widetilde X_{\hat0},\widetilde Y_{\hat0})\in(\widetilde D^N_1)
_{\hat0}\subset(\widetilde E_{10})_{\hat0}$.  Also, since $\widetilde D^T_{10}$ is totally geodesic by Lemma 2.8, we have $\nabla^{\widetilde M}
_{\widetilde X_{\hat0}}\widetilde Y\in(\widetilde D^T_{10})_{\hat0}\subset
(\widetilde E_{10})_{\hat0}$.  Therefore, we have $\widetilde h_{10}
(\widetilde X_{\hat0},\widetilde Y_{\hat0})=0$.  Thus we have 
$$\widetilde h_{10}((\widetilde E^T_{10})_{\hat0},
(\widetilde E^T_{10})_{\hat0})=0.\leqno{(2.6)}$$

Next we show $\widetilde h_{10}((\widetilde E^T_{10})_{\hat0},
(\widetilde D^N_1)_{\hat0})=0$.  Let $w\in(D^N_1)_{eK}$.  According to 
$(2.4)$, we suffices to show that $\widetilde h_{10}(\widehat X,\widehat w)
\,\,(X\in(D^T_{10})_{eK}),\,\,\widetilde h_{10}(\widehat{\eta},\widehat w)\,\,
(\eta\in\mathfrak c_{\mathfrak f}((D^N_2)_{eK}))$, 
$\widetilde h_{10}({\it l}^{co}_{Z,k},\widehat w)$ and 
$\widetilde h_{10}({\it l}^{si}_{Z,k},\widehat w)\,\,
(Z\in\mathfrak c_{\mathfrak g}((D^N_2)_{eK}),\,\,k\in{\bf N}\setminus\{0\})$ 
vanish.  According to Lemmas 2.2 and 2.3 in [Koi2], we have 
$\widetilde{\nabla}_{\widehat X}\overline w^L=(\nabla_X\overline w)_{\hat0}^L
-t[X,w],\,\,
\widetilde{\nabla}_{\widehat{\eta}}\overline w^L=-t[\eta,w],\,\,
\widetilde{\nabla}_{{\it l}^{co}_{Z,k}}
\overline w^L=-[\int_0^t{\it l}^{co}_{Z,k}(t)dt,\,w]$ and 
$\widetilde{\nabla}_{{\it l}^{si}_{Z,k}}
\overline w^L=-[\int_0^t{\it l}^{si}_{Z,k}(t)dt,\,w]$.  
Also, we can show 
$\nabla_X\overline w=-A_wX\in(D^T_{10})_{eK},\,\,[X,w]\in
\mathfrak c_{\mathfrak f}((D^N_2)_{eK})$, $[\eta,w]\in
\mathfrak c_{\mathfrak g}((D^N_2)_{eK})\cap((D^T_{10})_{eK}
\oplus(D^N_1)_{eK})$ and 
$[\int_0^t{\it l}^{co}_{Z,k}(t)dt,w],\,[\int_0^t{\it l}^{si}_{Z,k}(t)dt,w]
\in\mathfrak c_{\mathfrak g}((D^N_2)_{eK})\cap((D^T_{10})_{eK}\oplus
(D^N_1)_{eK})$ for each fixed 
$t\in[0,1]$.  Hence it follows from $(2.4)$ that 
$\widetilde{\nabla}_{\widehat X}\overline w^L ,\,
\widetilde{\nabla}_{\widehat{\eta}}\overline w^L$, 
$\widetilde{\nabla}_{{\it l}^{co}_{Z,k}}\overline w^L$ and 
$\widetilde{\nabla}_{{\it l}^{si}_{Z,k}}\overline w^L$ belong to 
$(\widetilde E_{10})_{\widehat 0}$.  That is, 
we have $\widetilde h_{10}(\widehat X,\widehat w)=\widetilde h_{10}
(\widehat{\eta},\widehat w)=\widetilde h_{10}({\it l}^{co}_{Z,k},\widehat w)
=\widetilde h_{10}({\it l}^{si}_{Z,k},\widehat w)=0$.  
Thus we have 
$$\widetilde h_{10}((\widetilde E^T_{10})_{\hat0},(\widetilde D^N_1)_{\hat0})
=0.\leqno{(2.7)}$$
Similarly, we can show $\widetilde h_{10}((\widetilde D^N_1)_{\hat0},
(\widetilde E^T_{10})_{\hat0})=0$, 
which together with $(2.5)\sim(2.7)$ and $(\widetilde E_{10})_{\hat0}
=(\widetilde E^T_{10})_{\hat0}\oplus(\widetilde D^N_1)_{\hat0}$ implies that 
$(\widetilde h_{10})_{\hat0}=0$, that is, $\hat 0$ is a geodesic point of 
$\widetilde E_{10}$.  This completes the proof of the totally geodesicness 
of $\widetilde E_{10}$.  Similarly, we can show that $\widetilde E_{20}$ and 
$\widetilde E_i$ ($i=1,2$) are totally geodesic.  
\hspace{7.5truecm}q.e.d.

\vspace{0.5truecm}

Let $\widetilde M_i(u):=\widetilde M\cap(u+V_i)$ and $(F_i)_u:=
T_u\widetilde M_i(u)$ ($u\in\widetilde M,\,i=1,2$).  

\vspace{0.5truecm}

{\bf Lemma 2.11.} {\sl The correspondence $F_i\,:\,u\,\mapsto\,(F_i)_u$ 
$(u\in \widetilde M)$ gives a totally geodesic distribution on $\widetilde M$ 
having $\widetilde M_i(u)$'s $(u\in\widetilde M)$ as integral manifolds, 
where $i=1,2$.}

\vspace{0.5truecm}

{\it Proof.} Fix $u_0\in\widetilde M$.  From (ii) of Lemma 2.6, 
it follows that 
$V_i=\overline{{\rm Span}\displaystyle{\left(\mathop{\cup}_{u\in 
L^{\widetilde D^T_i}_{u_0}}(\widetilde D^N_i)_u\right)}}$, 
where $L^{\widetilde  D^T_i}_{u_0}$ is the leaf of 
$\widetilde D^T_i$ through $u_0$.  On the other hand, it follows from Lemma 
2.10 that $(\widetilde D^N_i)_u$'s ($u\in L^{\widetilde D^T_i}_{u_0}$) 
are contained in 
$T_{u_0}L^{\widetilde D^T_i}_{u_0}\oplus(\widetilde D^N_i)_{u_0}$.  Hence we 
have $V_i\subset T_{u_0}
L^{\widetilde D^T_i}_{u_0}\oplus(\widetilde D^N_i)_{u_0}$ and hence 
$\widetilde M_i(u_0)\subset L^{\widetilde D^T_i}_{u_0}$.  It is clear that 
$\widetilde M_i(u_0)$ is totally geodesic in $L^{\widetilde D^T_i}_{u_0}$.  
Also, according to Lemma 2.8, $L^{\widetilde D^T_i}_{u_0}$ is totally 
geodesic in $\widetilde M$.  Hence $\widetilde M_i(u_0)$ is totally 
geodesic in $\widetilde M$.  This completes the proof.  
\hspace{0.1truecm}q.e.d.

\vspace{0.5truecm}

By using this lemma, we can show the following fact.  

\vspace{0.5truecm}

{\bf Lemma 2.12.} {\sl The submanifold $\widetilde M_i(u)$'s {\rm(}$u\in
\widetilde M${\rm)} are integral manifolds of $\widetilde D^T_i$ 
{\rm(}$i=1,2${\rm)}.}

\vspace{0.5truecm}

{\it Proof.} Let $\widetilde M'(u):=\widetilde M\cap(u+V')$ 
($u\in\widetilde M$).  
Since $V'=V_1\oplus V_2$ (orthogonal direct sum) by Lemma 2.9, 
we have $T_u\widetilde M'(u)=(F_1)_u\oplus(F_2)_u$ (orthogonal direct sum) 
for each $u\in\widetilde M$.  Also, it follows from Lemma 2.7 that 
$T_u\widetilde M'(u)=(\widetilde D^T_1)_u\oplus(\widetilde D^T_2)_u$ 
(orthogonal direct sum) for each $u\in\widetilde M$.  On the other hand, it 
follows from the proof of Lemma 2.11 that $(F_i)_u\subset(\widetilde D^T_i)_u$ 
($u\in\widetilde M,\,i=1,2$).  These facts imply $F_i=\widetilde D^T_i$ 
($i=1,2$).  Hence the statement of this lemma follows.  \hspace{1.6truecm}
q.e.d.

\vspace{0.5truecm}

\centerline{
\unitlength 0.1in
\begin{picture}( 51.3100, 20.2600)(  0.8000,-28.2600)
%
\special{pn 8}%
\special{pa 2140 1360}%
\special{pa 1244 2512}%
\special{fp}%
\special{pa 1244 2512}%
\special{pa 4444 2512}%
\special{fp}%
\special{pa 4444 2512}%
\special{pa 5212 1360}%
\special{fp}%
\special{pa 5212 1360}%
\special{pa 2140 1360}%
\special{fp}%
%
\special{pn 20}%
\special{sh 1}%
\special{ar 3676 2000 10 10 0  6.28318530717959E+0000}%
\special{sh 1}%
\special{ar 3676 1994 10 10 0  6.28318530717959E+0000}%
%
\special{pn 8}%
\special{pa 2140 1616}%
\special{pa 4488 2192}%
\special{fp}%
\special{pa 2678 1450}%
\special{pa 4316 2340}%
\special{fp}%
%
\special{pn 8}%
\special{pa 2012 2000}%
\special{pa 4578 2000}%
\special{fp}%
%
\special{pn 8}%
\special{ar 3676 2256 1070 768  3.2530750 4.2350640}%
%
\special{pn 8}%
\special{pa 2844 990}%
\special{pa 2844 1968}%
\special{fp}%
\special{pa 2844 2046}%
\special{pa 2844 2468}%
\special{fp}%
\special{pa 2844 2558}%
\special{pa 2844 2788}%
\special{fp}%
%
\special{pn 20}%
\special{sh 1}%
\special{ar 2844 1784 10 10 0  6.28318530717959E+0000}%
\special{sh 1}%
\special{ar 2844 1784 10 10 0  6.28318530717959E+0000}%
%
\special{pn 8}%
\special{pa 2748 1700}%
\special{pa 2652 1784}%
\special{pa 2978 1880}%
\special{pa 3062 1796}%
\special{pa 3062 1790}%
\special{pa 3062 1790}%
\special{pa 2748 1700}%
\special{fp}%
%
\special{pn 8}%
\special{pa 2984 1552}%
\special{pa 2876 1616}%
\special{pa 3150 1764}%
\special{pa 3266 1700}%
\special{pa 3266 1700}%
\special{pa 2984 1552}%
\special{fp}%
%
\special{pn 8}%
\special{pa 2562 1944}%
\special{pa 2492 2058}%
\special{pa 2824 2058}%
\special{pa 2876 1944}%
\special{pa 2876 1944}%
\special{pa 2562 1944}%
\special{fp}%
%
\special{pn 8}%
\special{pa 4066 2148}%
\special{pa 3958 2212}%
\special{pa 4232 2360}%
\special{pa 4348 2296}%
\special{pa 4348 2296}%
\special{pa 4066 2148}%
\special{fp}%
%
\special{pn 8}%
\special{pa 4168 2058}%
\special{pa 4072 2142}%
\special{pa 4398 2238}%
\special{pa 4482 2154}%
\special{pa 4482 2148}%
\special{pa 4482 2148}%
\special{pa 4168 2058}%
\special{fp}%
%
\special{pn 8}%
\special{pa 4252 1936}%
\special{pa 4182 2052}%
\special{pa 4514 2052}%
\special{pa 4566 1936}%
\special{pa 4566 1936}%
\special{pa 4252 1936}%
\special{fp}%
%
\special{pn 8}%
\special{pa 3420 2768}%
\special{pa 2742 2058}%
\special{da 0.070}%
\special{sh 1}%
\special{pa 2742 2058}%
\special{pa 2774 2120}%
\special{pa 2778 2098}%
\special{pa 2802 2092}%
\special{pa 2742 2058}%
\special{fp}%
\special{pa 3420 2768}%
\special{pa 2946 1872}%
\special{da 0.070}%
\special{sh 1}%
\special{pa 2946 1872}%
\special{pa 2960 1940}%
\special{pa 2972 1920}%
\special{pa 2996 1922}%
\special{pa 2946 1872}%
\special{fp}%
\special{pa 3420 2768}%
\special{pa 3150 1752}%
\special{da 0.070}%
\special{sh 1}%
\special{pa 3150 1752}%
\special{pa 3148 1822}%
\special{pa 3164 1804}%
\special{pa 3186 1810}%
\special{pa 3150 1752}%
\special{fp}%
\special{pa 3420 2768}%
\special{pa 4008 2244}%
\special{da 0.070}%
\special{sh 1}%
\special{pa 4008 2244}%
\special{pa 3946 2274}%
\special{pa 3968 2280}%
\special{pa 3972 2304}%
\special{pa 4008 2244}%
\special{fp}%
\special{pa 3420 2768}%
\special{pa 4258 2192}%
\special{da 0.070}%
\special{sh 1}%
\special{pa 4258 2192}%
\special{pa 4192 2214}%
\special{pa 4214 2222}%
\special{pa 4214 2246}%
\special{pa 4258 2192}%
\special{fp}%
\special{pa 3426 2768}%
\special{pa 4310 2052}%
\special{da 0.070}%
\special{sh 1}%
\special{pa 4310 2052}%
\special{pa 4246 2078}%
\special{pa 4268 2086}%
\special{pa 4270 2110}%
\special{pa 4310 2052}%
\special{fp}%
%
\special{pn 8}%
\special{pa 3062 1296}%
\special{pa 3150 1604}%
\special{da 0.070}%
\special{sh 1}%
\special{pa 3150 1604}%
\special{pa 3152 1534}%
\special{pa 3136 1554}%
\special{pa 3112 1546}%
\special{pa 3150 1604}%
\special{fp}%
%
\special{pn 8}%
\special{pa 2690 1624}%
\special{pa 2838 1776}%
\special{da 0.070}%
\special{sh 1}%
\special{pa 2838 1776}%
\special{pa 2806 1714}%
\special{pa 2800 1738}%
\special{pa 2776 1742}%
\special{pa 2838 1776}%
\special{fp}%
\special{pa 3804 1834}%
\special{pa 3670 1994}%
\special{da 0.070}%
\special{sh 1}%
\special{pa 3670 1994}%
\special{pa 3728 1956}%
\special{pa 3704 1954}%
\special{pa 3696 1930}%
\special{pa 3670 1994}%
\special{fp}%
\put(34.7000,-28.2600){\makebox(0,0)[rt]{$\widetilde E_i$}}%
\put(27.1500,-15.2700){\makebox(0,0)[rt]{$u$}}%
\put(32.2400,-11.3800){\makebox(0,0)[rt]{$\widetilde M_i(u)$}}%
\put(49.7400,-14.2400){\makebox(0,0)[rt]{$u+V_i$}}%
\put(43.1000,-8.0000){\makebox(0,0)[rt]{$u+V_0\oplus V_j$ ($j\in\{1,2\}\setminus\{i\}$)}}%
\put(41.4400,-16.5800){\makebox(0,0)[rt]{$\notin \widetilde U$}}%
\end{picture}%
\hspace{2.5truecm}
}

\vspace{0.5truecm}

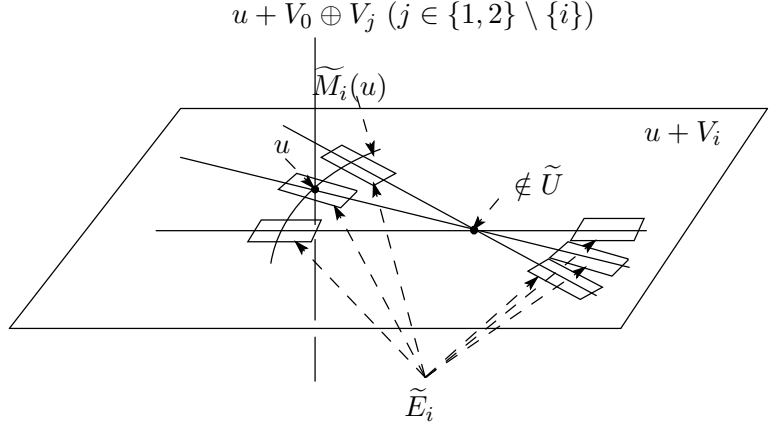
\centerline{{\bf Fig. 3.}}

\vspace{0.5truecm}

By using Lemma 2.11, we can show the following fact.  

\vspace{0.5truecm}

{\bf Lemma 2.13.} {\sl For any two points $u_1$ and $u_2$ of 
$\widetilde M'$, $\widetilde M_1(u_1)$ intersects with 
$\widetilde M_2(u_2)$.}

\vspace{0.5truecm}

{\it Proof.} Denote by $\mathfrak F_1$ the foliation on $\widetilde M'$ 
consisting of the integral manifolds of $F_1\vert_{\widetilde M'}$.  
Since $\mathfrak F_1$ is totally geodesic by Lemma 2.11 and the induced metric 
on each leaf of $\mathfrak F_1$ is complete, $F_2\vert_{\widetilde M'}$ is 
an Ehresmann connection for $\mathfrak F_1$ in the sense of Blumenthal-Hebda 
and hence the statement of this lemma follows (see [BH]).  
\hspace{4.7truecm}q.e.d.

\vspace{0.5truecm}

By using this lemma and imitating the proof of Corollary 3.11 of [HL], 
we can show the following fact.  

\vspace{0.5truecm}

{\bf Lemma 2.14.} {\sl For any $u_0\in\widetilde M_i(=\widetilde M_i(\hat 0))$, the translation map $f_{u_0}:V'\to V'$ defined by $f_{u_0}(u):=u+u_0$ 
$(u\in V')$ maps $\widetilde M_j(=\widetilde M_j(\hat0))$ 
isometrically onto $\widetilde M_j(u_0)$, where $(i,j)=(1,2)$ or $(2,1)$.}

\vspace{0.5truecm}

By using this lemma and imitating the proof of Corollary 3.12 of [HL], we 
can show the following fact.  

\vspace{0.5truecm}

{\bf Proposition 2.15.} {\sl We have $\widetilde M'=\widetilde M_1\times
\widetilde M_2\subset V_1\times V_2=V'$.}

\vspace{0.5truecm}

Define ideals $\mathfrak g'$ and $\mathfrak g_i$ ($i=1,2$) by 
$$\begin{array}{l}
\displaystyle{\mathfrak g':={\rm Span}\mathop{\cup}_{x^{\ast}\in M^{\ast}}
\{g_{0\ast}v(x^{\ast})_{\ast}^{-1}g_{0\ast}^{-1}\,\vert\,v\in 
T^{\perp}_{x^{\ast}}M^{\ast},\,g_0\in G\},}\\
\displaystyle{\mathfrak g_i:={\rm Span}\mathop{\cup}_{x^{\ast}\in M^{\ast}}
\{g_{0\ast}v(x^{\ast})_{\ast}^{-1}g_{0\ast}^{-1}\,\vert\,v\in((D^N_i)_
{\pi(x^{\ast})})_{x^{\ast}}^L,\,\,g_0\in G\}.}
\end{array}$$
Also, set $\mathfrak g_0:=\mathfrak g\ominus\mathfrak g'$, which is also an 
ideal of $\mathfrak g$.  Let $G'$ and $G_i$ ($i=0,1,2$) be the connected Lie 
subgroups of $G$ whose Lie algebras are $\mathfrak g'$ and $\mathfrak g_i$ 
($i=0,1,2$), respectively.  
Since $G/K$ is simply connected, we may assume that $G$ is simply connected.  
So we have $G=G'\times G_0$ and $G'=G_1\times G_2$.  By imitating the proof of 
Lemma 5.1 of [Koi4], we can show the following fact.  

\vspace{0.5truecm}

{\bf Lemma 2.16.} {\sl We have $V'\subset H^0([0,1],\mathfrak g')$ and 
$V_i\subset H^0([0,1],\mathfrak g_i)$ {\rm(}$i=1,2${\rm)}.}

\vspace{0.5truecm}

Also, by using Lemma 2.9 and imitating the proof of Lemma 3.7 of [E1], 
we can show the following fact.  

\vspace{0.5truecm}

{\bf Lemma 2.17.} {\sl We have $\mathfrak g_1\perp\mathfrak g_2$ and hence 
$H^0([0,1],\mathfrak g')=H^0([0,1],\mathfrak g_1)$\newline
$\oplus H^0([0,1],\mathfrak g_2)$ $($orthogonal direct sum$)$.}

\vspace{0.5truecm}

Let $V'_0:=H^0([0,1],\mathfrak g')\ominus V'$ and $V_{i,0}:=H^0([0,1],
\mathfrak g_i)\ominus V_i$ ($i=1,2$).  Clearly we have $V_0'=V_{1,0}\oplus 
V_{2,0}$.  Set $\widetilde M'_{H^0}:=\widetilde M\cap H^0([0,1],\mathfrak g')$ 
and $\widetilde M_{i,H^0}:=\widetilde M\cap H^0([0,1],\mathfrak g_i)$ 
($i=1,2$).  It follows from Proposition 2.1 that $\widetilde M'_{H^0}=
\widetilde M'\times V_0'$ and $\widetilde M_{i,H^0}=\widetilde M_i\times 
V_{i,0}$ ($i=1,2$).  
Furthermore, it follows from Proposition 2.15 that $\widetilde M=
\widetilde M_{1,H^0}\times\widetilde M_{2,H^0}\times 
H^0([0,1],\mathfrak g_0)$.  
It is clear that 
the parallel transport map $\phi$ for $G$ is decomposed as $\phi
=\phi_1\times\phi_2\times\phi_0$, where $\phi_i$ ($i=0,1,2$) is the 
parallel transport map for $G_i$.  Set $M^{\ast}_{i,H^0}:=\phi_i
(\widetilde M_{i,H^0})$ ($i=1,2$).  Clearly we have $M^{\ast}=
M^{\ast}_{1,H^0}\times M^{\ast}_{2,H^0}\times G_0\subset G_1\times G_2
\times G_0=G$.  Let $(\mathfrak g,\theta)$ be the orthogonal symmetric Lie 
algebra of $G/K$.  By imitating the discussion in Section 4 of [E1], we can 
show the following fact.  

\vspace{0.5truecm}

{\bf Lemma 2.18.} {\sl We have $\theta(\mathfrak g_i)=\mathfrak g_i$ 
{\rm(}$i=0,1,2${\rm)}.}

\vspace{0.5truecm}

Let $\mathfrak f_i:={\rm Fix}(\theta\vert_{\mathfrak g_i})$ and $K_i:=
\exp_{G_i}(\mathfrak f_i)$, where $i=0,1,2$.  Since $G/K$ is simply 
connected, we have $G/K=G_1/K_1\times G_2/K_2\times G_0/K_0$.  Denote by 
$\pi_i$ the natural projection of $G_i$ onto $G_i/K_i$ ($i=0,1,2$).  Let 
$M_{i,H^0}:=\pi_i(M^{\ast}_{i,H^0})$ ($i=1,2$).  
Now we prove Theorem A.  

\vspace{0.5truecm}

{\it Proof of Theorem A.}  Assume that the holonomy group of the section 
$\Sigma$ is reducible.  Then, under the above notations, we have 
$M=M_{1,H^0}\times M_{2,H^0}\times G_0/K_0\subset G_1/K_1\times G_2/K_2
\times G_0/K_0=G/K$.  
Let $\mathfrak t:=T^{\perp}_{eK}M$ and $\mathfrak t_i$ ($i=1,2$) be 
the normal space of $M_{i,H^0}$ in $G_i/K_i$.  Since $M$ is equifocal, 
$\mathfrak t$ is a Lie triple system.  Hence it follows that 
$\mathfrak t_i$ ($i=1,2$) are Lie triple systems.  This fact implies that 
$M_{i,H^0}$ ($i=1,2$) have Lie triple systematic normal bundle.  
On the other hand, it is clear that 
$M_{i,H^0}$ ($i=1,2$) satisfy the conditions (PF-i) and (PF-ii).  
Thus $M_{i,H^0}$ ($i=1,2$) is equifocal.  
The converse is trivial.  
\hspace{3.1truecm}q.e.d.

Next we shall prove Theorem B in terms of Theorem A.  

\vspace{0.5truecm}

{\it Proof of Theorem B.} Let $\Sigma$ be the section of $M$ through 
$x_0=g_0K\,(\in M)$ 
and $\pi_{\Sigma}:\widehat{\Sigma}\to\Sigma$ be the universal covering of 
$\Sigma$.  
Since $G/K$ is irreducible, it follows from Theorem A that 
the holonomy group of $\Sigma$ is irreducible, that is, $\widehat{\Sigma}$ 
is irreducible.  Since $\Sigma$ is totally geodesic in $G/K$, it is a 
symmetric space.  Hence $\widehat{\Sigma}$ is an irreducible simply connected 
symmetric space.  
On the other hand, according to Lemma 1A.4 of 
[PoTh1], $\Sigma$ and hence $\widehat{\Sigma}$ admit a totally geodesic 
hypersurface.  Hence, it follows from the result in [CN] that 
$\widehat{\Sigma}$ is isometric to a sphere, that is, $\Sigma$ is 
isometric to a sphere or a real projective space (of constant curvature).  
\hspace{8.9truecm}q.e.d.

\vspace{1truecm}

\centerline{{\bf References}}

\vspace{0.5truecm}

{\small 
\noindent
[A] M. M. Alexandrino, Singular riemannian foliations with sections, 
Illinois J. 

Math. {\bf 48} (2004) 1163-1182.

\noindent
[AG] M. M. Alexandrino and C. Gorodski, Singular Riemannian foliations 
with 

sections, transnormal maps and basic forms, arXiv:math.GT/0608069.

\noindent
[AT] M. M. Alexandrino and D. T$\ddot o$ben, Singular Riemannian foliations on 

simply connected spaces, Differential Geom. Appl. {\bf 24} (2006) 383-397.

\noindent
[B] L. Biliotti, Coisotropic and polar actions 
on compact irreducible Hermitian 

symmetric spaces, Trans. Amer. Math. Soc. {\bf 358} (2006) 3003-3022.

\noindent
[BH] R. A. Blumenthal and J. J. Hebda, 
Complementary distributions which 

preserves the leaf geometry and applications to totally geodesic foliations, 

Quat. J. Math. {\bf 35} (1984) 383-392.

\noindent
[BV] J. Berndt and L. Vanhecke, 
Curvature-adapted submanifolds, 
Nihonkai 

Math. J. {\bf 3} (1992) 177-185.

\noindent
[BCO] J. Berndt, S. Console and C. Olmos, Submanifolds and holonomy, 
Re-

search Notes in Mathematics 434, CHAPMAN $\&$ HALL/CRC Press, Boca 

Raton, London, New York 
Washington, 2003.

\noindent
[BS] R. Bott and H. Samelson, Applications of the theory of Morse to 
symmetric 

spaces, Amer. J. Math. {\bf 80} (1958) 964-1029. Correction in Amer. J. Math. 

{\bf 83} (1961) 207-208.  

\noindent
[Ch] U. Christ, 
Homogeneity of equifocal submanifolds, J. Differential Geome-

try {\bf 62} (2002) 1-15.

\noindent
[CKT] Y.W. Choe, U.H. Ki and R. Takagi, 
Compact minimal generic subman-

ifolds with parallel normal section in a complex projective space, Osaka J. 

Math. {\bf 37} (2000) 489-499. 

\noindent
[CN] B.Y. Chen and T. Nagano, 
Totally geodesic submanifolds of symmetric 

spaces II, Duke Math. J. {\bf 45} (1978) 405-425.

\noindent
[Co] L. Conlon, 
Remarks on commuting involutions, 
Proc. Amer. Math. Soc. 

{\bf 22} (1969) 255-257.

\noindent
[D] J. Dadok, Polar coordinates induced by actions of compact Lie 
groups, 

Trans. Amer. Math. Soc. {\bf 228} (1985) 125-137.

\noindent
[E1] H. Ewert, 
A splitting theorem for equifocal 
submanifolds in simply con-

nected compact symmetric spaces, Proc. Amer. Math. Soc. {\bf 126} 
(1998) 

2443-2452.

\noindent
[E2] H. Ewert, Equifocal submanifolds in Riemannian symmetric spaces, 
Doc-

toral thesis.

\noindent
[G] C. Gorodski, 
Polar actions on compact symmetric spaces which admit a 

totally geodesic principal orbit, Geometriae Dedicata {\bf 103} (2004) 
193-204.  

\noindent
[H] S. Helgason, 
Differential geometry, Lie groups and symmetric spaces, Aca-

demic Press, New York, 1978.

\noindent
[HL] E. Heintze and X. Liu, 
A splitting theorem for isoparametric submanifolds 

in Hilbert space, J. Differential Geom. {\bf 45} (1997) 319-335.

\noindent
[HLO] E, Heintze, X. Liu and C. Olmos, Isoparametric submanifolds and a 

Chevalley-type restriction theorem, arXiv:math.DG/0004028.

\noindent
[HPTT] E. Heintze, R.S. Palais, C.L. Terng and G. Thorbergsson, 
Hyperpolar 

actions on symmetric spaces, Geometry, topology and physics for Raoul 

Bott (ed. S. T. Yau), Conf. Proc. Lecture 
Notes Geom. Topology {\bf 4}, 

Internat. Press, Cambridge, MA, 1995 pp214-245.

\noindent
[Koi1] N. Koike, Tubes of non-constant radius in symmetric spaces, Kyushu J. 

Math. {\bf 56} (2002) 267-291.

\noindent
[Koi2] N. Koike, On proper Fredholm submanifolds in a Hilbert space arising 

from submanifolds in a symmetric space, Japan. J. Math. {\bf 28} (2002) 61-80.

\noindent
[Koi3] N. Koike, Actions of Hermann type and proper complex equifocal 
sub-

manifolds, Osaka J. Math. {\bf 42} (2005) 599-611.

\noindent
[Koi4] N. Koike, A splitting theorem for proper complex equifocal 
submanifolds, 

Tohoku Math. J. {\bf 58} (2006) 393-417.

\noindent
[Kol1] A. Kollross, A Classification of hyperpolar and cohomogeneity one 
ac-

tions, Trans. Amer. Math. Soc. {\bf 354} (2001) 571-612.

\noindent
[Kol2] A. Kollross, Polar actions on symmetric spaces, J. Differential Geomery 

{\bf 77} (2007) 425-482.

\noindent
[KW] A. Kor$\grave a$nyi and J.A. Wolf, Realization of Hermitian symmetric 
spaces as 

generalized half-planes, Ann. of Math. {\bf 81} (1965) 265-288.

\noindent
[LT] A. Lytchak and G. Thorbergsson, Variationally complete actions on 
non-

negatively curved manifolds, Illinois J. Math. (to appear).

\noindent
[M] P. Molino, Riemannian foliations, Progress in Mathematics vol. 73, 
Birkh$\ddot a$-

user Boston 1988. 

\noindent
[Pa] R.S. Palais, 
Morse theory on Hilbert manifolds, Topology {\bf 2} (1963) 299-340.

\noindent
[PaTe] R.S. Palais and C.L. Terng, Critical point theory and submanifold 
ge-

ometry, Lecture Notes in Math. {\bf 1353}, Springer, Berlin, 1988.

\noindent
[PoTh1] F. Podest$\acute a$ and G. Thorbergsson, 
Polar actions on rank-one symmetric 

spaces, J. Differential Geometry {\bf 53} (1999) 131-175.

\noindent
[PoTh2] F. Podest$\acute a$ and G. Thorbergsson, 
Polar and coisotropic actions on 

K$\ddot a$hler manifolds, 
Trans. Amer. Math. Soc. {\bf 354} (2002) 1759-1781.

\noindent
[Ta] R. Takagi, Real hypersurfaces in a complex projective space with 
constant 

principal curvatures, J. Math. Soc. Japan {\bf 27} (1975) 43-53. 

\noindent
[TaTa] R. Takagi and T. Takahashi, On the principal curvatures of 
homoge-

neous hypersurfaces in a sphere, Differential Geometry, in honor of 

K. Yano, Kinokuniya, Tokyo, 1972, 469-481.

\noindent
[TeTh] C.L. Terng and G. Thorbergsson, 
Submanifold geometry in symmetric 

spaces, J. Differential Geometry {\bf 42} (1995) 665-718.

\noindent
[Th] G. Thorbergsson, Isoparametric foliations and their buildings, 
Ann. of 

Math. {\bf 133} (1991) 429-446.

\noindent
[T$\ddot o$] D. T$\ddot o$ben, Parallel focal structure and singular 
Riemannian foliations, 

Trans. Amer. Math. Soc. {\bf 358} (2006) 1677-1704.



\noindent
[Ts] K. Tsukada, Totally geodesic hypersurfaces of naturally reductive 
homo-

geneous spaces, Osaka J. Math. {\bf 33} (1996) 697-707.

\vspace{0.5truecm}

{\small 
\rightline{Department of Mathematics, Faculty of Science}
\rightline{Tokyo University of Science, 26 Wakamiya-cho}
\rightline{Shinjuku-ku, Tokyo 162-8601 Japan}
\rightline{(koike@ma.kagu.tus.ac.jp)}
}
\end{document}